\documentclass[twoside,11pt]{amsart}

\usepackage{latexsym,amsmath,amsfonts}

\usepackage{amsmath}
\usepackage{latexsym}
\usepackage{amsthm}
\usepackage{amssymb}
\usepackage{graphicx}

\renewcommand{\theequation}{\arabic{section}.\arabic{equation}}

\newtheorem{theorem}{Theorem}[section]
\newtheorem{remark}{Remark}[section]
\newtheorem{coro}{Corollary}[section]
\newtheorem{lemma}{Lemma}[section]
\newtheorem{prop}{Proposition}[section]
\numberwithin{equation}{section}

\numberwithin{equation}{section}

\begin{document}

\title{\bf Global existence and asymptotic decay of solutions
to the Non-isentropic Euler-Maxwell system}

\author{Yuehong Feng$^{1}$, Shu Wang$^{1}$ and Shuichi
Kawashima$^{2}$}

\date{}

\maketitle \markboth{Y. H. Feng, S. Wang and S. Kawashima } { Global
existence and asymptotic decay of solutions to the Non-isentropic
Euler-Maxwell system}

\vspace{-3mm}

\begin{center}
{\small $^1$College of Applied Sciences, Beijing University of Technology, Beijing 100124, China \\[2mm]
$^2$Faculty of Mathematics, Kyushu University, Fukuoka, 812-8581, Japan \\[2mm]
Email~: fengyuehong001@163.com, \hspace{1mm}wangshu@bjut.edu.cn,
\hspace{1mm} kawashim@math.kyushu-u.ac.jp}
\end{center}


\vspace{1cm}

\begin{center}
\begin{minipage}{14cm}
    {\bf Abstract.} {\small The non-isentropic compressible Euler-Maxwell system is
investigated in $R^3$ in the present paper, and the $L^q$ time decay
rate for the global smooth solution is established. It is shown that
the density and temperature of electron converge to the equilibrium
states at the same rate $(1+t)^{-\frac{11}{4}}$ in $L^q$ norm.}
\end{minipage}
\end{center}

\vspace{7mm}


\vspace{5mm}

\noindent {\bf AMS Subject Classification (2000)~:} 35A01, 35L45,
35L60, 35Q35


\vspace{3mm}


\section{Introduction and main results}

We investigate the asymptotic behavior of globally smooth solutions
for the (rescaled) non-isentropic Euler-Maxwell systems, which takes
the following form(see Appendix A and
\cite{Chen84,Dink05,Jero03,Jero05,YW11}):
\begin{equation}
\label{1.1}
\left\{\begin{aligned}
&\partial_t n + \nabla\cdot(n  u )=0, \\
&\partial_t u + ( u\cdot\nabla) u +\nabla \theta+\theta\nabla \ln n
+ u = -(E+   u \times B)
, \\
&\partial_t \theta +  u\cdot\nabla\theta
+\frac{2}{3}\theta\nabla\cdot u
= \frac{1}{3}|u|^2-(\theta-1) , \\
&  \partial_t E-\nabla\times B=   n u ,
 \\
&  \partial_t B+\nabla\times E=0,\\
&     \nabla\cdot E=1-n ,\quad \nabla\cdot B=0,\quad
(t,x)\in(0,\infty)\times\mathbb{R}^3.
\end{aligned} \right.
\end{equation}
Here, $n,u,\theta$ denote the scaled macroscopic density, mean
velocity vector and temperature of the electrons and $E,B$ the
scaled electric field and magnetic field. They are functions of a
three-dimensional position vector $x\in\mathbb{R}^3$ and of the time
$t>0$. The fields $E$ and $B$ are coupled to the particles through
the Maxwell equations and act on the particles via the Lorentz force
$E+ u \times B$.

In this paper, we are interested in the asymptotics and global
existence of smooth solutions of system (\ref{1.1}) with the initial
conditions~:
\begin{equation}
\label{1.2}
    (n, u, \theta, E,B)|_{t=0} = (n_0, u_0, \theta_0, E_0,B_0),\quad x\in
    \mathbb{R}^3,
\end{equation}
which satisfies the compatible condition
\begin{equation}
   \label{1.3}
    \nabla\cdot E_0 = 1 - n_0, \quad \nabla\cdot B_0 = 0, \quad x \in \;
    \mathbb{R}^3.
\end{equation}

The Euler-Maxwell system (\ref{1.1}) is a symmetrizable hyperbolic
system for $n,\theta >0$. It is known that the Cauchy problem
(\ref{1.1})-(\ref{1.2}) has a local smooth solution when the initial
data are smooth. In a simplified one dimensional isentropic
Euler-Maxwell system, the global existence of entropy solutions has
been given in \cite{{CJW00}} by the compensated compactness method.
For the three dimensional isentropic Euler-Maxwell system, the
existence of global smooth solutions with small amplitude to the
Cauchy problem in the whole space and to the periodic problem in the
torus is established by Peng et al in \cite{PWG11} and Ueda et al in
\cite{UKW10} respectively, and the decay rate of the smooth solution
when $t$ goes to infinity is obtained by Duan in \cite{Duan11} and
Ueda et al in \cite{UK11}. For asymptotic limits with small
parameters, see \cite{PW08a,PW08b} and references therein. Recently,
Yang et al in \cite{YW11} consider the diffusive relaxation limit of
the three dimensional non-isentropic Euler-Maxwell system.

In this paper, we consider large time asymptotics and global
existence of the smooth solutions to the non-isentropic
Euler-Maxwell system. Our main goal here is to establish the global
existence of smooth solutions around a constant state $(1,0,1,0,0)$,
which is a equilibrium solution of system (\ref{1.1}), and the decay
rate of the global smooth solutions in time for the system
(\ref{1.1}). Our main results read as follows:
\begin{theorem}\label{thm1.1}
Let $s\geq 4$ be an integer, and \eqref{1.3} hold. There exist
$\delta_0>0$ and a constant $C_0$ such that if $$\left\|
{[n_{0}-1,u_{0},\theta_0-1,E_0,B_0] } \right\|_s\leq \delta_0,$$
then, the initial problem \eqref{1.1}- \eqref{1.2} has a unique
global solution $[n (t,x),u (t,x)$,
 $\theta(t,x), E(t,x)$, $B(t,x)]$ with $$[n -1,u,\theta-1, E, B]\in C^1\big([0,T);H^{s-1}(\mathbb{R}^3)\big) \cap
C\big([0,T);H^s(\mathbb{R}^3)\big)$$ and
 $$
\mathop {\sup }\limits_{t \geqslant 0}\left\|
{[n(t)-1,u(t),\theta(t)-1,E(t),B(t)] } \right\|_s\leq C_0\left\|
{[n_{0}-1,u_{0},\theta_0-1,E_0,B_0] } \right\|_s.
$$ Furthermore, there exist $\delta_1>0$ and a constant $C_1$ such that if
$$\left\| {[n_{0}-1,u_{0},\theta_0-1,E_0,B_0] }
\right\|_{13}+\left\| {[u_{0},E_0,B_0] } \right\|_{L^1}\leq
\delta_1,$$ then $[n (t,x),u (t,x)$,
 $\theta(t,x), E(t,x), B(t,x)]$ satisfies
\begin{equation}\notag
\begin{split}
& \left\|\left[
{n(t)-1, \theta(t)-1 }\right] \right\|_{L^q}\leq C_1 (1+t)^{-\frac{11}{4}}, \\
&\left\|
{ E(t)  } \right\|_{L^q}\leq C_1 (1+t)^{-2+\frac{3}{2q}}, \\
& \left\|
{u(t),B}(t)  \right\|_{L^q}\leq C_1 (1+t)^{-\frac{3}{2}+\frac{3}{2q}}, \\
\end{split}
\end{equation}with  $t\geq 0$ and $2\leq q \leq\infty$.
\end{theorem}
The proof of Theorem \ref{thm1.1} is based on the careful energy
methods and the Fourier multiplier technique. This is divided into
three key steps: The first key step is to establish the global a
priori high order Sobolev's energy estimates in time by using the
careful energy methods and the skey-symmetric dissipative structure
(see \cite{UKW10}) of Euler-Maxwell system, which concludes the
global existence results in Theorem \ref{thm1.1}. Due to the
complexity of non-isentropic case caused by the coupled energy
equations, the technique of symmetrizer is used here to remove this
difficulty when we deal with the Euler part of the Euler-Maxwell
system so as to obtain the energy estimates for the the density, the
velocity and the temperature while the special symmetric structure
of the Maxwell system is used when we obtain the energy estimates
for the electric-magnetic fields for the Maxwell part of the
Euler-Maxwell system, which is different from the method used by
\cite{Duan11, UKW10} to deal with the isentropic case of
compressible Euler-Maxwell system. The second key step is to obtain
the $L^p-L^q$ time decay rate of the linearized operator for the
non-isentropic Euler-Maxwell system by using the Fourier technique,
used first by Kawashima in \cite{SK84} and extended then to the
other problems, see \cite{DH95, Duan11, UK11} and the references
therein. Here we first apply energy method in the Fourier space to
obtain the basic $L^\infty$ estimates for the Fourier transform of
the solution. Then we solve the dissipative linear wave system of
three order by the Fourier technique. We find that the 'error'
functions $\rho$, $\Theta$, $\nabla\cdot u$ satisfy the same
dissipative linear wave equation which is of order three and is
different from that of the isentropic Euler-Maxwell system. We also
find an interesting phenomenon that the estimate on $B$ which is
depending on the temperature $\Theta$ is different from that of the
isentropic Euler-Maxwell equation. Since the linear system involved
here is of order three, it is complex to obtain the time decay rate
of the linear system by the Fourier analysis. Fortunately, we can
obtain the elaborate spectrum structure of the eigenvalue equations
of this linear system of order three, which yield to the desired
time decay rate for the linearized system. The third step of the
proof is to obtain the time decay rate in Theorem \ref{thm1.1} by
combining the previous two steps and apply the energy estimate
technique to the nonlinear problem satisfied by the error functions,
whose solutions can be represented by the solution-semigroup
operator for the linearized problem by using the Duhamel Principle.
This concludes the time asymptotic stability results in Theorem
\ref{thm1.1}

For the later use in this paper, we give some notations. For any
integer $s\geq0$, $H^s$, $\dot{H}^s$ denote the Sobolev space
$H^s({{\mathbb R}^3}) $ and the s-order homogeneous Sobolev space,
respectively. Set $L^2=H^0$. The norm of $H^s$ is denoted by
$\left\| {\cdot} \right\|_s$ with $\left\| {\cdot} \right\|=\left\|
{\cdot} \right\|_0$. $f \sim g$ denotes $\gamma f\leq g\leq
\frac{1}{\gamma}f$ for a constant $0<\gamma<1$.
$\langle\cdot,\cdot\rangle$ denotes the inner product in
$L^2({\mathbb R}^3)$, i.e.
$$ \langle f,g\rangle=
\int_{{\mathbb R}^3} {f(x)g(x)dx,} \quad  f=f(x),~g=g(x)\in
L^2({\mathbb R}^3).$$ We set $\partial^\alpha$
$=\partial_{x_1}^{\alpha_1}
\partial_{x_2}^{\alpha_2}\partial_{x_3}^{\alpha_3}$ $=\partial_{1}^{\alpha_1}
\partial_{2}^{\alpha_2}\partial_{3}^{\alpha_3}$ for a multi-index $\alpha
=[\alpha_1,\alpha_2,\alpha_3]$,  For an integrable function
$f:{\mathbb R}^3\rightarrow{\mathbb R}$, its Fourier transform is
defined by
$$ \hat{f}(k)=\int_{{\mathbb R}^3}e^{-ix\cdot k}f(x)dx,\ \ x\cdot k:=
\sum\limits_{j=1}^{3} {x_j k_j},\ \ k\in {\mathbb R}^3,
$$
where $i=\sqrt{-1}\in \mathbb C$ is the imaginary unit. For two
complex numbers or vectors $a$ and $b$, $(a|b):=a\cdot\overline{b}$
denotes the dot product of $a$ with the complex conjugate of $b$.
 We also use
$\mathcal {R}(a)$  to denote the real part of $a$.

This paper is organized as follows. In Section 2, the transformation
of system \eqref{1.1} is presented. In Section 3, we obtain the
existence and uniqueness of global solutions. In Section 4, we study
the linearized homogeneous equations to get the $L^p-L^q$ decay
property and the explicit representation of solutions. Lastly, in
Section 5, we investigate the decay rates of solutions of the
nonlinear system \eqref{2.2} and complete the proof of Theorem
\ref{thm1.1}.
\section{Transformation of system \eqref{1.1}}
Suppose $[n (t,x), u (t,x),\theta(t,x),E(t,x),B(t,x)]$ to be a
smooth solution to the initial problem of the non-isentropic
Euler-Maxwell equations \eqref{1.1} with initial data \eqref{1.2}
which satisfies \eqref{1.3}. We introduce the transformation
\begin{equation}\label{reformulation}
n (t,x)=1+\rho (t,x), \theta(t,x)=1+\Theta (t,x).
\end{equation}
Then, the system \eqref{1.1} becomes
\begin{equation}
\label{2.2}
\left\{\begin{aligned}
&\partial_t \rho + \nabla\cdot((1+\rho)  u )=0, \\
&\partial_t u + ( u\cdot\nabla) u +\nabla
\Theta+\frac{1+\Theta}{1+\rho}\nabla \rho  = -(E+   u \times B)-u
, \\
&\partial_t \Theta +  u\cdot\nabla\Theta
+\frac{2}{3}(1+\Theta)\nabla\cdot u
= \frac{1}{3}|u|^2-\Theta , \\
&  \partial_t E-\nabla\times B=   (1+\rho) u ,
 \\
&  \partial_t B+\nabla\times E=0,\\
& \nabla\cdot E=-\rho ,\quad \nabla\cdot B=0,\quad
(t,x)\in(0,\infty)\times\mathbb{R}^3,
\end{aligned} \right.
\end{equation}
with initial data
\begin{equation}\label{2.3}
U|_{t=0}=[\rho,u,\Theta,E,B]|_{t=0}=U_0:=[\rho_{ 0},u_{ 0},\Theta_{
0},E_0,B_0],\ x \in{\mathbb R}^3,
\end{equation}
which satisfies the compatible condition
\begin{equation}\label{2.4}
    \nabla\cdot E_0 = -\rho_0, \quad \nabla\cdot B_0 = 0, \quad x \in  \;
    \mathbb{R}^3,
\end{equation}
with $\rho_0=n_0-1.$ In the following, we set $s\geq4$. Besides, for
$U=[\rho,u,\Theta,E,B]$, we use $\mathcal {E}_s(U(t))$, $\mathcal
{E}_s^h(U(t))$, $\mathcal {D}_s(U(t))$, $\mathcal {D}_s^h(U(t))$ to
denote the energy functional, the high-order energy functional, the
dissipation rate and the high-order dissipation rate, respectively.
Here,
\begin{equation}\label{2.5}
\mathcal {E}_s(U(t))\sim \left\| {[\rho,u,\Theta,E,B]} \right\|_s^2,
\end{equation}
\begin{equation}\label{2.6}
\mathcal {E}_s^h(U(t))\sim \left\| {\nabla[\rho,u,\Theta,E,B]}
\right\|_{s-1}^2,
\end{equation}
\begin{equation}\label{2.7}
\mathcal {D}_s(U(t))\sim \left\| {[\rho,u,\Theta]} \right\|_s^2+
\left\| {\nabla[E,B]} \right\|_{s-2}^2+\left\| {E} \right\|^2,
\end{equation}
and
\begin{equation}\label{2.8}
\mathcal {D}_s^h(U(t))\sim \left\| {\nabla[\rho,u,\Theta]}
\right\|_{s-1}^2+ \left\| {\nabla[E,B]} \right\|_{s-2}^2.
\end{equation}
Then, for the initial problem \eqref{2.2}-\eqref{2.3}, we obtain
\begin{prop}\label{prop2.1}
Assume that $U_0=[\rho_{0},u_{0},\Theta_{0},E_0,B_0]$ satisfies
\eqref{2.4}. Then, there exist $ \mathcal {E}_s(\cdot)$ and
$\mathcal {D}_s(\cdot)$ such that the following holds true. If \ $
\mathcal {E}_s(U_0)>0$ is small enough, then, for any $t \geq0$, the
initial problem \eqref{2.2}-\eqref{2.3} has a unique nonzero global
solution $U =[\rho,u,\Theta,E,B]$ which satisfies
\begin{equation}\label{2.9}
U\in C^1\big([0,T);H^{s-1}(\mathbb{R}^3)\big) \cap
C\big([0,T);H^s(\mathbb{R}^3)\big),
\end{equation}
and
\begin{equation}\label{2.10}
\mathcal {E}_s(U(t))+\gamma\int_0^t\mathcal {D}_s(U(y))dy\leq
\mathcal {E}_s(U_0).
\end{equation}

\end{prop}
Furthermore, we investigate the time decay rates of solutions in
Proposition \ref{prop2.1} under some extra conditions on the given
initial data $U_0=$$[\rho_{0}$, $u_{0}$, $\Theta_{0}$, $E_0$,
$B_0]$. For this purpose, we define $\epsilon_s(U_0)$ as
\begin{equation}\label{2.11}
\epsilon_s(U_0)=\left\|{U_0}\right\|_s+\left\|{[u_{0},E_0,B_0]}\right\|_{L^1}
\end{equation}
for $s\geq4.$ Then, we have
\begin{prop}\label{prop2.2}
Assume that \eqref{2.4} holds for given initial data
$U_0=[\rho_{0},u_{0},\Theta_{0},E_0,B_0]$. If $ \epsilon_{s+2}(U_0)$
is small enough, then $U =[\rho,u,\Theta,E ,B ]$ satisfies
\begin{equation}\label{2.12}
\left\|{U(t)}\right\|_s\leq C
\epsilon_{s+2}(U_0)(1+t)^{-\frac{3}{4}}
\end{equation}
for $t\geq0$. Moreover, if $\epsilon_{s+6}(U_0)$ is small enough,
then $U$ $=$[$\rho,u$, $\Theta$, $E$, $B$ $]$ also satisfies
\begin{equation}\label{2.13}
\left\|{\nabla U(t)}\right\|_{s-1}\leq C
 \epsilon_{s+6}(U_0)(1+t)^{-\frac{5}{4}}
\end{equation}
for $t \geq0$.
\end{prop}

\begin{prop}\label{prop2.3}
Assume that \eqref{2.4} holds for given initial data
$U_0=$$[\rho_{0}$, $u_{0}$, $\Theta_{0}$, $E_0$, $B_0]$. If\ $
\epsilon_{13}(U_0)$ is small enough, then $U =[\rho,u,\Theta,E ,B ]$
satisfies
\begin{equation}\label{2.14}
\left\|{[\rho(t),\Theta(t)]}\right\|_{L^q}\leq C
(1+t)^{-\frac{11}{4}},
\end{equation}
\begin{equation}\label{2.15}
\left\|{[u(t),E(t)]}\right\|_{L^q}\leq C (1+t)^{-2+\frac{3}{2q}},
\end{equation}
and
\begin{equation}\label{2.16}
\left\|{B(t)}\right\|_{L^q}\leq C (1+t)^{-\frac{3}{2}+\frac{3}{2q}}
\end{equation}
for $2\leq q\leq \infty$ and $t\geq0$.
\end{prop}
Lastly, one can obtain Theorem \ref{thm1.1} from Proposition
\ref{prop2.1} and Proposition \ref{prop2.3}. Therefore, we prove the
three previous Propositions in the rest of this paper.
%

\section{Global solutions for system \eqref{2.2}.} We review
Moser-type calcluls inequalities in Sobolev spaces and
 the local existence of smooth solutions for symmetrizable hyperbolic
equations, which will be used in the proof of our main theorem.

\begin{lemma}
\label{L3.1}
    (Moser-type calculus inequalities, see \cite{KM81,Ma84})
Let $s \geq 1$ be an integer. Suppose $u \in H^s(\mathbb{R}^3)$,
$\nabla u \in L^\infty(\mathbb{R}^3)$ and $v \in
H^{s-1}(\mathbb{R}^3) \cap L^\infty(\mathbb{R}^3)$. Then for all
multi-index $\alpha$ with $|\alpha|\leq s$, we have
$\partial^{\alpha} (uv)-u\partial^{\alpha}  v\in L^2(\mathbb{R}^3)$
and
$$  \|\partial^\alpha (uv)-u\partial^{\alpha}  v\| \leq C_s\big(\|\nabla u\|_{L^\infty} \|D^{s-1} v\|
+ \|D^s u\| \|v\|_{L^\infty}\big),  $$ where
$$  \|D^s u\| = \sum\limits_{|\alpha| = s} \|\partial^{\alpha}  u\|.    $$
Moreover, if $s > \frac{5}{2}$, then the embedding
$H^{s-1}(\mathbb{R}^3) \hookrightarrow L^\infty(\mathbb{R}^3)$ is
continuous and we have
$$  \|\partial^\alpha (uv)-u\partial^{\alpha}  v\| \leq C_s \|u\|_s \|v\|_{s-1},
\quad \|uv\|_{s-1} \leq C_s \|u\|_{s-1} \|v\|_{s-1}.    $$
\end{lemma}

\begin{lemma}
\label{L3.2}
    (Local existence of smooth solutions, see \cite{Ka75,Ma84})
Let $s > \frac{5}{2}$ and $(\rho_{ 0},u_{ 0},\Theta_{ 0},E_0,B_0)
\in H^s(\mathbb{R}^3)$. Then there exist $T> 0$ and a unique smooth
solution $(n, u, \theta, E, B)$ to the Cauchy problem
(\ref{1.1})-(\ref{1.2}) satisfying $(\rho, u, \Theta, E, B) \in
C^1\big([0,T);H^{s-1}(\mathbb{R}^3)\big) \cap
C\big([0,T);H^s(\mathbb{R}^3)\big)$.
\end{lemma}
\noindent 3.1. \textbf{Preliminary results.} In this subsection, we
will prove the following a priori estimates
\begin{theorem}\label{thm3.1}
  Assume that
$U =[\rho,u,\Theta,E ,B ]\in
C^1\big([0,T);H^{s-1}(\mathbb{R}^3)\big) \cap
C\big([0,T);H^s(\mathbb{R}^3)\big) $ is smooth for $T>0$ with
\begin{equation}\label{3.1}
\mathop {\sup }\limits_{0\leq t\leq T
}\left\|{U(t)}\right\|_s\leq\delta
\end{equation}
for $\delta\leq \delta_0$ with $\delta_0$ sufficiently small and
suppose $U$ to be the solution of the system \eqref{2.2} for
$t\in(0,T)$. Then, there exist $\mathcal {E}_s(\cdot)$ and
$\mathcal{D}_s(\cdot)$ such that
\begin{equation}\label{3.2}
\frac{d}{dt}\mathcal {E}_s(U(t))+\gamma \mathcal {D}_s(U(t))\leq
C[\mathcal {E}_s(U(t))^\frac{1}{2}+\mathcal {E}_s(U(t))]\mathcal
{D}_s(U(t))
\end{equation}
for $0\leq t\leq T.$
\end{theorem}

\noindent \emph{Proof.} We will use five steps to finish the proof.
The step 1 is to estimate the Euler part and the Maxwell part of the
Euler-Maxwell system, respectively. Steps 2-4 is to obtain the
dissipative estimates for $\rho$, $E$ and $B$
by using the skew-symmetric structure of the Euler-Maxwell system.\\
\emph{Step 1.} It holds that
\begin{equation}\label{3.3}
\frac{d}{dt}\sum_{|\alpha|\leq s} [\langle A^{I
}_0(W_I)\partial^\alpha W_{I},\partial^\alpha
    W_{I}\rangle+\|\partial^\alpha W_{II}\|^2
    ]+\|u\|_s^2+\frac{1}{3}\|\Theta\|_s^2\leq
    C\|W\|_s\|W_I\|_s^2,\ \ \
\end{equation}
where
$$  W_{I} = \left(\begin{array}{c}
    \rho \\[2mm] u \\[2mm] \Theta
\end{array} \right),
 \quad W_{II}=\left(\begin{array}{c}
    E \\[2mm] B
\end{array} \right), \quad W =\left(\begin{array}{c}
    W_{I }\\[2mm] W_{II}
\end{array} \right), \quad A_0^I({W_I}) = \left( {\begin{array}{*{20}{c}}
   {\frac{{1 + \Theta }}{{1 + \rho}}} & 0 & 0  \\
   0 & {(1 + \rho){I_3}} & 0  \\
   0 & 0 & {\frac{3}{2}\frac{{1 + \rho}}{{1 + \Theta }}}  \\
\end{array}} \right).
$$
In fact, Set
\[A_j^I({W_I}) = \left( {\begin{array}{*{20}{c}}
   {{u_j}} & {(1 + \rho)e_j^T} & 0  \\
   {\frac{{1 + \Theta }}{{1 + \rho}}{e_j}} & {{u_j}{I_3}} & {{e_j}}  \\
   0 & {\frac{2}{3}(1 + \Theta )e_j^T} & {{u_j}}  \\
\end{array}} \right)\begin{array}{*{20}{c}}
   {} & {j = 1,2,3,}  \\
\end{array}\]
$$  K_{I}(W) = \left(\begin{array}{c}
    0 \\[2mm] -(E + u \times B)\\[2mm]0
\end{array} \right), \quad
    K_2(W) = \left(\begin{array}{c}
    0 \\[2mm] -u\\[2mm] \frac{1}{3}|u|^2-\Theta
\end{array} \right).
$$
Then the first three equations of \eqref{2.2} for $W_{I }$ can be
rewritten under the form
\begin{equation}
\label{3.4}
    \partial_t W_{I } +\sum\limits_{j=1}^3 A^{I }_j(W_{I })\partial_{x_j}W_{I } = K_{I }(W) +  K_2(W).
\end{equation}
It is clear that system \eqref{3.4} for $W_{I }$ is symmetrizable
hyperbolic when $1+\rho,1+\Theta > 0$. More precisely, since we
consider small solutions defined in a time interval $[0,T)$ with $T
> 0$, \eqref{3.1} implies that $\|\left[\rho,\Theta\right]\|_{L^\infty((0,T) \times \mathbb{R}^3)}
\leq C_s\|\left[\rho,\Theta\right] \|_s \leq C_s \delta \leq
\frac{1}{2}$. Then $\frac{1}{2} \leq 1+\rho,1+\Theta  \leq
\frac{3}{2}$. It follows that $A_0^I({W_I})$ is symmetric positively
definite and $\tilde A^{I }_j(W_{I })=A^{I }_0(W_I)A^{I }_j(W_{I})$
are symmetric for all $1 \leq j \leq 3$. This choice of $A^{I
}_0(W_I)$ will simplify energy estimates.

For $|\alpha|\leq s$, differentiating equations \eqref{3.4} with
respect to $x$ and multiplying the resulting equations by the
symmetrizer matrix $A^{I}_0(W_I)$, one has
\begin{equation}
\label{3.5}
\begin{split}
 A^{I }_0(W_I) \partial_t\partial ^\alpha W_{I} +
\sum\limits_{j=1}^3
    A^{I }_0(W_I)A_j^{I }(W_{I })\partial_{x_j} \partial ^\alpha W_{I }
 =  A^{I }_0(W_I)\partial ^\alpha \Big(K_{I }(W)+ K_2(W)\Big)
+J_\alpha,
\end{split}\end{equation}
where $J_\alpha$ is defined by
$$  J_\alpha=-\sum_{j=1}^3 A^{I }_0(W_I)\big[\partial ^\alpha
    \big(A_j^{I }(W_{I })\partial_{x_j} W_{I }\big)-A_j^{I }(W_{I })\partial ^\alpha (\partial_{x_j}W_{I })\big].   $$
Applying Lemma \ref{L3.1} to $J_\alpha$, we get
\begin{equation}
\label{3.6}
\begin{split}
    \|J_\alpha\| & \leq  C \big(\|\nabla A^{I}_j(W_{I})\|_{L^\infty} \|\partial_{x_j} W_{I}\|_{s-1}
+ \|D^s A^{I}_j(W_{I})\| \|\partial_{x_j} W_{I}\|_{L^\infty}\big) \\
    & \leq C\big(\|W_{I}\|_s + \|W_{I}\|_s\big)\|W_{I}\|_s  \leq C\|W_{I}\|^2_s.
\end{split}\end{equation} Taking the inner product of equations \eqref{3.5}
with $\partial^\alpha W_{I}$ and using the fact that the matrix
$\tilde A^{I}_j(W_{I})$ is symmetric, we have
\begin{equation}
\label{3.7}
\begin{split}    \frac{d}{dt} \langle A^{I }_0(W_I)\partial ^\alpha W_{I},\partial ^\alpha
    W_{I}\rangle
 = & 2 \langle J_\alpha,\partial ^\alpha W_{I}\rangle + \langle
\mathop{\mathrm{div}}
 A^{I}(W_{I})\partial ^\alpha W_{I},\partial ^\alpha W_{I}\rangle \\
& + 2 \langle A^{I}_0(W_I)\partial ^\alpha W_{I},\partial^\alpha
K_{I}(W)+ \partial^\alpha  K_2(W)\rangle,
\end{split}\end{equation}
where
\begin{eqnarray}
    \mathop{\mathrm{div}} A^{I}(W_{I})  = \partial_t A^{I}_0(W_I)
    + \sum_{j=1}^3 \partial_{x_j} \tilde{A}_j^{I}(W_{I})
 = \partial_\rho A^{I }_0(W_{I}) \partial_t \rho+\partial_\Theta A^{I }_0(W_{I}) \partial_t \Theta + \sum_{j=1}^3
(\tilde{A}^{I}_j)'(W_{I}) \partial_{x_j} W_{I}. \nonumber
\end{eqnarray}
Using the first equation in \eqref{2.2} and Lemma \ref{L3.1}, we
have
\begin{equation}
\label{3.8}
    \|\partial_t \rho\|_{L^\infty} \leq  C \|\partial_t \rho\|_{s-1}
= C \big\|\nabla\cdot \big((1+\rho)u\big)\big\|_{s-1}\leq  C
\big(1+\|\rho\|_s\big) \|u\|_s.
\end{equation}
Then
\begin{equation}
\label{3.9}
    \|\mathop{\mathrm{div}} A^{I}(W_{I})\|_{L^\infty} \leq C(1+\|W_{I}\|_s)\|W_{I}\|_s.
\end{equation}

Now let us estimate each term on the right hand side of \eqref{3.7}.
For the first two terms, by Cauchy-Schwarz inequality and using
estimates \eqref{3.6} and \eqref{3.9}, we have
\begin{equation}
\label{3.10}
\begin{split}
    \langle J_\alpha,\partial^\alpha  W_{I}\rangle
+ \langle\mathop{\mathrm{div}} A^{I}(W_{I})\partial^\alpha
W_{I},\partial^\alpha W_{I}\rangle
 \leq
C(1+\|W_{I}\|_s)\|W_{I}\|^3_s \leq C \|W\|_s\|W_{I}\|^2_s.
\end{split}
\end{equation}
For the third term of the right hand side of \eqref{3.7}, by using
the definition of $A^{I}_0(W_I)$, $K_{I}(W)$ and $K_2(W),$  we
obtain
\begin{equation}
\label{3.11}
\begin{split}& 2 \langle A^{I }_0(W_I)\partial^\alpha W_{I},\partial^\alpha
K_{I}(W)
+ \partial^\alpha K_2(W)\rangle \\
& =- 2 \langle(1+\rho)\partial^\alpha  u,\partial^\alpha u\rangle -
2\langle\partial^\alpha  u,\partial^\alpha  E\rangle -
2\langle\partial^\alpha
 u,\partial^\alpha  (u \times B)\rangle  - 2\langle \rho \partial^\alpha  u , \partial^\alpha  E
\rangle \\
 & \quad -2\langle \rho \partial^\alpha  u , \partial^\alpha  (u
\times B)\rangle
 +\langle \frac{1+\rho}{2(1+\Theta)}\partial^\alpha
\Theta,\partial^\alpha (|u|^2)\rangle
-\langle \frac{3(1+\rho)}{2(1+\Theta)}\partial^\alpha \Theta,\partial^\alpha \Theta\rangle \\
& \leq  - \|\partial^\alpha  u\|^2 -\frac{1}{3}\|\partial^\alpha
\Theta\|^2- 2\langle\partial^\alpha  u,\partial^\alpha  E\rangle +
C\|W\|_s \|W_{I}\|_s^2.
\end{split}
\end{equation}
Next we write the system for $W_{II}$ as~:
$$  \partial_t W_{II} +\sum\limits_{j=1}^3 A^{II}_j\partial_{x_j}W_{II} = \big((1+\rho )u  ,0\big)^T,   $$
where

$$  A^{II}_j \!=\! \left(\begin{array}{cc}
    0 & L_j \\[2mm]
    L^T_j & 0
\end{array} \right), \; \; L_1 \!=\! \left( \begin{array}{ccc}
    0 & 0 & 0 \\
    0 & 0 & 1 \\
    0 & -1 & 0 \\
\end{array} \right), \; \; L_2 \!=\! \left( \begin{array}{ccc}
    0 & 0 & -1 \\
    0 & 0 & 0 \\
    1 & 0 & 0 \\
\end{array} \right), \; \; L_3 \!=\! \left( \begin{array}{ccc}
    0 & 1 & 0 \\
    -1 & 0 & 0 \\
    0 & 0 & 0 \\
\end{array} \right).
$$
For $|\alpha| \leq s$, differentiating the fourth and fifth
equations of (\ref{2.2}) with respect to $x$, we get
\begin{equation}
\label{3.12}
 \left\{\begin{split}
&\partial_t\partial^\alpha  E-\nabla\times\partial^\alpha  B =\partial^\alpha  \big[(1+\rho)u\big], \\
&\partial_t\partial^\alpha  B +\nabla\times\partial^\alpha  E=0.
\end{split}\right.
\end{equation}
By the vector analysis formula
$$  \nabla\cdot(f\times g)=(\nabla\times f)\cdot g-(\nabla\times g)\cdot f,  $$
one term appearing in Sobolev energy estimates vanishes, i.e.
$$
\int_{\mathbb{R}^3}\big(-\nabla\times\partial^\alpha  B\cdot
\partial^\alpha  E+\nabla\times\partial^\alpha
E\cdot\partial^\alpha B\big)dx =
\int_{\mathbb{R}^3}\nabla\cdot\big(\partial^\alpha
E\times\partial^\alpha B\big)dx = 0.
$$
Hence, the standard energy estimate for \eqref{3.12} together with
Lemma \ref{L3.1} yields
\begin{equation}
\label{3.13}
\begin{split}
    \frac{d}{dt}\|\partial ^\alpha W_{II}\|^2   =2\langle\partial^\alpha  E,\partial^\alpha  u\rangle
    +2\langle\partial^\alpha  E,\partial^\alpha  (\rho u)\rangle
     \leq 2\langle\partial^\alpha  E,\partial^\alpha
     u\rangle+C\|W\|_s
\|W_{I}\|_s^2.
\end{split}
\end{equation}
 Adding \eqref{3.7}, \eqref{3.10}, \eqref{3.11} and \eqref{3.13}, summing up for
 all $|\alpha| \leq s$, we get (\ref{3.3}).

\noindent \emph{Step 2.} It holds that
\begin{equation}\label{3.14}
  \begin{split}  \frac{d}{dt}\sum_{\alpha\leq s-1}\langle \frac{1}{2(1+\rho)}\partial ^\alpha\rho   -\partial
   ^\alpha\nabla\cdot u,
  \partial ^\alpha \rho\rangle+
    \gamma\|\rho\|_s^2  \leq C ( \|W\|_s  \|W_I\|_s^2+  \|[u,\Theta]\|_s^2 ).
        \end{split}
\end{equation}
In fact, for $|\alpha|\le s-1$. Differentiating the second equation
of \eqref{2.2} with respect to $x$ and taking the inner product of
the resulting equation with $\partial ^\alpha\nabla \rho$, we have
\begin{equation}
\label{3.15}
\begin{split}&\langle \frac{1+\Theta}{1+\rho}\partial ^\alpha\nabla \rho,
\partial ^\alpha\nabla \rho\rangle
+\langle\partial ^\alpha E,\partial^\alpha\nabla \rho\rangle  \\
& = -\langle\partial^\alpha(\frac{1+\Theta}{1+\rho}\nabla \rho) -
\frac{1+\Theta}{1+\rho}\partial^\alpha\nabla
\rho,\partial^\alpha\nabla \rho\rangle
-\langle\partial^\alpha\partial_tu ,\partial^\alpha\nabla \rho \rangle  \\
&\quad-\langle\partial^\alpha(u \cdot\nabla u +\nabla\Theta+ u\times
B),
\partial ^\alpha\nabla \rho\rangle - \langle\partial ^\alpha
u,\partial ^\alpha\nabla \rho\rangle.
\end{split}
\end{equation}
Let us estimate each term in \eqref{3.15}. First, noting that
$\frac{1}{2}\leq 1+\rho,1+\Theta \leq \frac{3}{2}$, we have
\begin{equation}
\label{3.16}
\begin{split} \langle \frac{1+\Theta}{1+\rho}\partial ^\alpha\nabla \rho,
\partial ^\alpha\nabla \rho\rangle
+\langle\partial ^\alpha E,\partial ^\alpha\nabla \rho\rangle =&
\langle \frac{1+\Theta}{1+\rho}\partial ^\alpha\nabla \rho,
\partial ^\alpha\nabla \rho\rangle+\langle\partial ^\alpha \rho, \partial ^\alpha \rho\rangle \\
  \geq & C^{-1} \big(\|\partial ^\alpha \nabla \rho\|^2 +
\|\partial ^\alpha \rho\|^2\big).
\end{split}
\end{equation}
By Lemma \ref{L3.1}, we obtain
\begin{equation}
\notag\begin{split} \|\partial^\alpha(\frac{1+\Theta}{1+\rho}\nabla
\rho) - \frac{1+\Theta}{1+\rho}\partial ^\alpha\nabla \rho\|
 \leq C\|W_I\|^2_s.
\end{split}
\end{equation}
Then,
\begin{equation}
\label{3.17}
    -\langle\partial^\alpha(\frac{1+\Theta}{1+\rho}\nabla \rho) -
\frac{1+\Theta}{1+\rho}\partial^\alpha\nabla
\rho,\partial^\alpha\nabla \rho\rangle \leq C\|W\|_s \|W_{I }\|^2_s.
\end{equation}
Obviously,
\begin{equation}
\notag\begin{split}
    -\langle\partial ^\alpha\partial_tu ,\partial ^\alpha\nabla \rho \rangle=
     \frac d{dt}\langle\partial ^\alpha \nabla\cdot u ,
    \partial ^\alpha
\rho \rangle -\langle\partial ^\alpha \nabla\cdot u ,\partial
^\alpha
\partial_t \rho \rangle.
\end{split}
\end{equation}
Then, from (\ref{3.8}) we have
\begin{equation}
\notag\begin{split}
    \big|\langle\partial ^\alpha \nabla\cdot u ,\partial^\alpha \partial_t \rho \rangle\big|
 \leq  \|\partial ^\alpha \nabla\cdot u \| \|\partial ^\alpha \partial_t \rho \|
\leq  \|u \|^2_s + \|W\|_s \|W_{I }\|_s^2.
\end{split}
\end{equation}
Hence,
\begin{equation}
\label{3.18}
    -\langle\partial ^\alpha\partial_tu ,\partial ^\alpha\nabla \rho \rangle
\leq \frac d{dt}\langle\partial ^\alpha \nabla\cdot u ,\partial
^\alpha \rho \rangle + \|u \|^2_s + \|W\|_s \|W_{I }\|_s^2.
\end{equation}
From (\ref{3.1}), one has
\begin{equation}
\label{3.19}
\begin{split}&  \big|\langle\partial ^\alpha(u \cdot\nabla u +\nabla\Theta+
u\times
B), \partial ^\alpha\nabla \rho\rangle\big|  \\
& \leq  \big(\|\partial ^\alpha(u \cdot\nabla u )\| + \|\partial
^\alpha (u \times B)\|
 +\|\partial ^\alpha\nabla\Theta\| \big) \|\partial ^\alpha \nabla \rho \|  \\
& \leq  C\big(\|u \|_{s-1} \|u \|_s + \|u \|_{s-1} \|B\|_{s-1}
+ \|u \|_{s-1}+\|\Theta\|_s \big) \|\rho \|_s \\
& \leq  C \|W\|_s \|W_{I }\|_s^2 + \varepsilon \|\rho \|^2_s +
C_\varepsilon \big(\|u\|^2_s+\|\Theta\|^2_s\big).
\end{split}
\end{equation}

Now we establish the uniform estimates for the last term on the
right hand side of (\ref{3.15}). We get
\begin{equation}
\label{3.20}
\begin{split}    -\langle\partial ^\alpha u,\partial ^\alpha\nabla \rho\rangle
 =&
 -\langle\partial ^\alpha(\frac {\partial_t
\rho}{1+\rho})-\frac{\partial ^\alpha \partial_t
\rho}{1+\rho},\partial ^\alpha \rho\rangle - \langle\frac {u\cdot
\partial ^\alpha \nabla \rho}{1+\rho},\partial ^\alpha \rho\rangle\\
&- \langle\frac {\partial_t\partial ^\alpha \rho}{1+\rho},\partial
^\alpha \rho\rangle
 -  \langle\partial ^\alpha(\frac{u\cdot\nabla
\rho}{1+\rho}) -\frac{u\cdot\partial ^\alpha \nabla
\rho}{1+\rho},\partial ^\alpha \rho\rangle,
\end{split}
\end{equation}
with
\begin{equation}
\label{3.21}
\begin{split}- \langle\frac {\partial_t\partial ^\alpha \rho }{1+\rho },\partial ^\alpha \rho \rangle
&=  - \frac{1}{2} \frac d{dt}\langle\frac 1{1+\rho }\partial ^\alpha
\rho ,\partial ^\alpha \rho \rangle
+ \frac{1}{2} \langle\partial_t(\frac{1}{1+\rho })\partial ^\alpha \rho ,\partial ^\alpha \rho \rangle  \\
&=  - \frac{1}{2} \frac d{dt}\langle\frac 1{1+\rho }\partial ^\alpha
\rho ,\partial ^\alpha \rho \rangle + \frac{1}{2}
\langle\frac{1}{(1+\rho )^2} \partial_t \rho  \partial ^\alpha \rho
,\partial ^\alpha \rho \rangle.
\end{split}
\end{equation}

Obviously,
\begin{equation}
\label{3.22}
    \Big|\langle\frac {u \cdot\nabla\partial ^\alpha \rho }{1+\rho },\partial ^\alpha \rho \rangle\Big|
 \leq  C \|\rho \|_{s-1} \|u \|_s \|\rho \|_s  \leq  C
\|W\|_s \|W_{I } \|^2_s.
\end{equation}
Using (\ref{3.8}) and \eqref{3.1}, we have
\begin{equation}
\label{3.23}
    \Big|\langle\frac{1}{(1+\rho )^2} \partial_t \rho  \partial^\alpha \rho ,\partial ^\alpha \rho \rangle \Big|
 \leq  C \|\partial_t \rho \|_{L^\infty} \|\partial ^\alpha \rho \|^2  \leq  C
\|W\|_s \|W_{I } \|^2_s.
\end{equation}
By Lemma \ref{L3.1} and the continuous embedding
$H^{s-1}(\mathbb{R}^3) \hookrightarrow L^\infty(\mathbb{R}^3)$, one
has
\begin{equation}
\notag\begin{split}
    \Big\|\partial ^\alpha(\frac{\partial_t \rho }{1+\rho })-\frac {\partial ^\alpha \partial_t \rho }
    {1+\rho }\Big\|
& \leq  C \big(\big\|\nabla \big(\frac{1}{1+\rho
}\big)\big\|_{L^\infty} \|\partial_t \rho \|_{s-2}
+ \big\|D^{s-1}\big(\frac{1}{1+\rho }\big)\big\| \|\partial_t \rho \|_{L^\infty})  \\
& \leq  C \big(\|\nabla \rho \|_{s-1} \|\partial_t \rho \|_{s-2} +
\|\rho \|_{s-1} \|\partial_t \rho
\|_{s-1}\big)  \\
& \leq  C \|\rho \|_s \|\partial_t \rho \|_{s-1}  \\
& \leq  C \big(1+\|\rho \|_s\big) \|u \|_s \|\rho \|_s.
\end{split}
\end{equation}
Then, since \eqref{3.1}, we have
\begin{equation}
\label{3.24}
    \Big|\langle\partial^\alpha(\frac {\partial_t \rho }{1+\rho })-\frac{\partial^\alpha \partial_t \rho }{1+\rho },
\partial^\alpha \rho \rangle\Big| \leq C\|W\|_s\|W_{I } \|^2_s.
\end{equation}
Similarly,
\begin{equation}
\notag\begin{split}
    \big\|\partial ^\alpha(\frac {u \cdot\nabla \rho }{1+\rho })
    -\frac{u \cdot \partial ^\alpha \nabla \rho }{1+\rho }\big\|
  \leq C \big(\big\|\nabla \big(\frac{u }{1+\rho
}\big)\big\|_{L^\infty} \|\nabla \rho \|_{s-2} +
\big\|D^{s-1}\big(\frac{u }{1+\rho }\big)\big\| \|\nabla \rho
\|_{L^\infty}\big)   \leq C \|W_{I }\|_s \|\rho \|_s.
\end{split}
\end{equation}
Therefore,
\begin{equation}
\label{3.25}
    \Big|\langle\partial^\alpha(\frac{u \cdot\nabla \rho }{1+\rho })
    -\frac{u \cdot \partial^\alpha \nabla \rho }{1+\rho },\partial ^\alpha \rho \rangle\Big|
 \leq C \|W\|_s \|W_{I }\|^2_s.
\end{equation}

Thus, combining (\ref{3.15})-(\ref{3.25}), we get
\begin{equation}\notag
\begin{split}
      \frac d{dt}\langle\frac 1{1+\rho }\partial ^\alpha \rho -\partial^\alpha  \nabla\cdot u,
\partial ^\alpha \rho \rangle  +  C^{-1} \big(\|\partial ^\alpha \nabla \rho  \|^2 + \|\partial ^\alpha
    \rho \|^2\big)    \leq  C\|W\|_s \|W_{I }\|^2_s + \varepsilon \|\rho \|^2_s +
C_\varepsilon\big( \|u \|^2_s+\|\Theta \|^2_s \big).
\end{split}
\end{equation}
Summing up this inequality for all $|\alpha| \leq s-1$ and taking
$\varepsilon > 0$ small enough, we obtain \eqref{3.14}.

\noindent \emph{Step 3.} It holds that
\begin{equation}\label{3.26}
  \begin{split}  \frac{d}{dt}\sum_{|\alpha|\leq s-1}\langle \partial ^\alpha u,\partial ^\alpha E\rangle  +
    \gamma\|E\|_{s-1}^2  \leq & C  \|[E,B]\|_s(\|u\|_s^2+  \|E\|_{s-2}^2)\\
     & +C\|[\rho,u,\Theta]\|_s^2+C \|u\|_s\|\nabla B\|_{s-2}.
        \end{split}
\end{equation}
In fact, for $|\alpha|\leq s-1$, applying $\partial^\alpha$ to the
second equation of \eqref{2.2}, multiplying it by $\partial^\alpha
E$, taking integration in $x$ and then using the fourth equation of
\eqref{2.2} implies
\begin{equation}\notag
\begin{split}
  &\frac{d}
{{dt}}\left\langle {{\partial ^\alpha }u,{\partial ^\alpha }E}
\right\rangle   + {\left\| {{\partial ^\alpha }E} \right\|^2} \\& =
- \left\langle {{\partial ^\alpha }(u \cdot \nabla u),{\partial
^\alpha }E} \right\rangle  - \left\langle {{\partial ^\alpha
}(\frac{{1 + \Theta }}{{1 + \rho }}\nabla \rho ),{\partial ^\alpha
}E} \right\rangle
    - \left\langle {\nabla {\partial ^\alpha }\Theta ,{\partial ^\alpha
    }E}
    \right\rangle - \left\langle {{\partial ^\alpha }u,{\partial ^\alpha }E}
     \right\rangle\\
    &\quad - \left\langle {{\partial ^\alpha }(u \times B),{\partial ^\alpha }E} \right\rangle
  + \left\langle {{\partial ^\alpha }u,\nabla  \times {\partial ^\alpha }B}
  \right\rangle + {\left\| {{\partial ^\alpha }u} \right\|^2} + \left\langle {{\partial ^\alpha }(\rho u),
    {\partial ^\alpha }u} \right\rangle.
\end{split}
\end{equation}
Furthermore, using Cauchy-Schwarz inequality, one has
\begin{equation}\notag
\begin{split}
 & \frac{d}
{{dt}}\left\langle {{\partial ^\alpha }u,{\partial ^\alpha }E}
\right\rangle
 + \gamma {\left\| {{\partial ^\alpha }E} \right\|^2}
 \\& =  C\left({\left\| E \right\|_{s - 1}}
 \left\| u \right\|_s^2+ \left\| {[\rho ,u,\Theta ]} \right\|_s^2  + {\left\| u \right\|_s}{\left\| {\nabla B}
    \right\|_{s - 2}}\right)
  +C{\left\| B \right\|_s}\left( {\left\| u \right\|_s^2 + \left\|
 E \right\|_{s - 1}^2} \right) \\
\end{split}
\end{equation}
Therefore, taking summation of the previous estimate over
$|\alpha|\leq s-1$, one has \eqref{3.26}.

\noindent \emph{Step 4.} It holds that
\begin{equation}\label{3.27}
\frac{d} {{dt}}\sum\limits_{\left| \alpha  \right| \leqslant s - 2}
{\left\langle { - \nabla  \times {\partial ^\alpha }E,{\partial
^\alpha }B} \right\rangle }  + \gamma \left\| {\nabla B} \right\|_{s
- 2}^2 \leqslant C\left\| {[u,E]} \right\|_{s - 1}^2 + C\left\| \rho
\right\|_{s}^2\left\| {\nabla u} \right\|_{{s - 1}}^2.\end{equation}

In fact, for $|\alpha|\leq s - 2$, applying $\partial^\alpha$ to the
fourth equation of \eqref{2.2}, multiplying it by
$-\partial^\alpha\nabla\times B $, taking integration in $x$ and
then using the fifth equation of \eqref{2.2} gives
\begin{equation}\notag
\begin{split}
   \frac{d}
{{dt}}\left\langle { - \nabla  \times {\partial ^\alpha }E,{\partial
^\alpha }B} \right\rangle
 + {\left\| {\nabla  \times {\partial ^\alpha }B} \right\|^2}  = {\left\| {\nabla  \times {\partial ^\alpha }E} \right\|^2} - \left\langle {{\partial ^\alpha }u,\nabla
 \times {\partial ^\alpha }B} \right\rangle  - \left\langle {{\partial ^\alpha }(\rho u),\nabla  \times
 {\partial ^\alpha }B} \right\rangle.
\end{split}
\end{equation}
Furthermore, using Cauchy-Schwarz inequality and taking summation
over $|\alpha|\leq s - 2$, one has \eqref{3.27}. Where we also used
$$\left\|{\partial^\alpha\partial_i B}\right\|= \left\|{\partial_i\triangle^{-1}\nabla\times(\nabla\times
\partial^\alpha B)}\right\|
\leq C\left\|{\nabla\times\partial^\alpha B}\right\|$$ for each
$1\leq i\leq 3$, due to the fact that $\nabla\cdot B=0$ and
$\partial_i\triangle^{-1}\nabla$ is bounded from $L^p$ to $L^p$ for
$1<p<\infty$, see \cite{Stein}.\\
\emph{Step 5.} Next, based four previous steps, we will prove
\eqref{3.2}. We define the energy functional as
\begin{equation}\notag
 \begin{split}
 \mathcal{E}_s(U(t))=& \sum_{|\alpha|\leq s} [\langle A^{I
}_0(W_I)\partial^\alpha W_{I},\partial^\alpha
    W_{I}\rangle+\|\partial^\alpha W_{II}\|^2
    ]\\
    & +\mathcal {K}_1 \sum\limits_{\left| \alpha  \right| \leqslant s -
1} \langle \frac{1}{2(1+\rho)}\partial ^\alpha\rho-\partial
^\alpha\nabla\cdot u,
  \partial ^\alpha \rho\rangle
+\mathcal {K}_2 \sum\limits_{\left| \alpha  \right| \leqslant s -
1}\langle \partial ^\alpha u,\partial ^\alpha E\rangle \\& +\mathcal
{K}_3 \sum\limits_{\left| \alpha  \right| \leqslant s - 2}
{\left\langle { - \nabla  \times {\partial ^\alpha }E,{\partial
^\alpha }B} \right\rangle },
  \end{split}
\end{equation}
where, constants $0<\mathcal {K}_3\ll\mathcal {K}_2\ll\mathcal
{K}_1\ll 1$ is to be chosen later. Notice the fact that $A^{I
}_0(W_I) $ is positively and that as soon as $0<\mathcal {K}_j
\ll1_{(j=1,2,3)}$ is sufficiently small, then $\mathcal
{E}_s(U(t))\sim ||U||^2_s$ holds true. Furthermore, the sum of
\eqref{3.3}, \eqref{3.14}$\times\mathcal {K}_1$,
\eqref{3.26}$\times\mathcal {K}_2$ and \eqref{3.27}$\times\mathcal
{K}_3$ gives that there exits $0<\gamma<1$ such that
\begin{equation}\notag
\begin{split}
\frac{d}{dt}&\mathcal {E}_s(U(t))+\|u\|_s^2
+\frac{1}{3}\|\Theta\|_s^2 + \gamma\mathcal {K}_1\|\rho\|_s^2+
\gamma\mathcal {K}_2\|E\|_{s-1}^2+ \gamma\mathcal {K}_3 \left\|
{\nabla B}
\right\|_{s - 2}^2 \\
&\leq
    C[\mathcal {E}_s(U(t))^\frac{1}{2}+\mathcal {E}_s(U(t))]\mathcal
{D}_s(U(t))+C \mathcal {K}_1 \|[u,\Theta]\|_s^2+C\mathcal {K}_2\|[\rho,u,\Theta]\|_s^2\\
&\quad  +C \mathcal {K}_3 \left\| {[u,E]} \right\|_{s - 1}^2+C
\mathcal {K}_2 \|u\|_s\|\nabla B\|_{s-2}\\
 &\leq
    C[\mathcal {E}_s(U(t))^\frac{1}{2}+\mathcal {E}_s(U(t))]\mathcal
{D}_s(U(t))+C \mathcal {K}_1 \|[u,\Theta]\|_s^2+C\mathcal {K}_2\|[\rho,u,\Theta]\|_s^2\\
&\quad  +C \mathcal {K}_3 \left\| {[u,E]} \right\|_{s - 1}^2+
\frac{1}{2}C\left(\mathcal {K}_2^\frac{1}{2}\left\|{u
}\right\|_s^2+\mathcal {K}_2^\frac{3}{2}\left\|{\nabla
B}\right\|_{s-2}^2\right),
\end{split}
\end{equation}
where the Cauchy-Schwarz inequality were used. By letting
$0<\mathcal {K}_3\ll\mathcal {K}_2\ll\mathcal {K}_1\ll 1$ be
sufficiently small with $\mathcal {K}_2^{\frac{3}{2}}\ll\mathcal
{K}_3$, we obtain that there exists $\gamma>0$, $C>0$ such that
\eqref{3.2} also holds true. Now, we have finished the proof of
Theorem \ref{thm3.1}. \hfill $\Box$ \vspace{0.2cm}

\noindent 3.2. \textbf{Proof of Proposition 2.1.} { } The global
existence of smooth solutions follows from the standard argument by
using the local existence result given in Lemma \ref{L3.2}, the a
priori estimate \eqref{3.2} given in Theorem \ref{thm3.1} and the
continuous extension argument, see \cite{{Ni78}}. The proof of the
Proposition \ref{prop2.1} is finished. \hfill $\Box$
%
\section{Linearized homogeneous equations of system \eqref{2.2}}
In this section, so as to obtain the time-decay rates of solutions
to the nonlinear system \eqref{2.2} in the next section, we consider
the following initial problem on the linearized homogeneous
equations corresponding to system \eqref{2.2}:
\begin{equation}
\label{4.1}
\left\{\begin{aligned}
&\partial_t \rho + \nabla\cdot u=0, \\
&\partial_t u  +\nabla \rho +\nabla \Theta+E+ u   =0, \\
&\partial_t \Theta+\frac{2}{3}\nabla\cdot u +\Theta= 0 , \\
&  \partial_t E-\nabla\times B-u =0, \\
&  \partial_t B+\nabla\times E=0,\\
& \nabla\cdot E=-\rho ,\quad \nabla\cdot B=0,\quad t > 0, x \in
\mathbb{R}^3,
\end{aligned} \right.
\end{equation}
with initial data
\begin{equation}\label{4.2}
U|_{t=0}=U_0:=[\rho_{ 0},u_{ 0},\Theta_{ 0},E_0,B_0],\ x \in{\mathbb
R}^3,
\end{equation}
which satisfies the compatible condition
\begin{equation}\label{4.3}
    \nabla\cdot E_0 = -\rho_0, \quad \nabla\cdot B_0 = 0, \quad x \in  \;
    \mathbb{R}^3.
\end{equation}
In this section, we usually use $U=[\rho, u, \Theta, E, B]$ to
denote the solution of the linearized homogeneous system
\eqref{4.1}.\vspace{0.3cm}

\noindent 4.1. \textbf{Pointwise time-frequency estimate.}
 In this subsection, we utilize the energy method to the initial
 problem \eqref{4.1}-\eqref{4.3} in the Fourier space to present that there is a time-frequency Lyapunov functional
  which is equivalent to $|\hat{U}(t,k)|^2$ and  furthermore its dissipation rate could be
  represented by itself. The main result of this subsection is presented in the following.
\begin{theorem}\label{thm4.1}
Suppose $U(t,x),(t,x)\in(0,\infty)\times\mathbb{R}^3 $ to be a
 solution of the linearized homogeneous system \eqref{4.1}. There exists a
time-frequency Lyapunov functional $ \mathcal {E}(\hat{U}(t,k))$
with
\begin{equation}\label{4.4}
\mathcal {E}(\hat{U})\sim
|\hat{U}|^2:=|\hat{\rho}|^2+|\hat{u}|^2+|\hat{\Theta}|^2+|\hat{E}|^2+|\hat{B}|^2
\end{equation}
which satisfies that there exists a constant $\gamma>0$ such that
\begin{equation}\label{4.5}
\frac{d}{dt}\mathcal {E}(\hat{U}(t,k))+\frac{\gamma
|k|^2}{(1+|k|^2)^2}\mathcal {E}(\hat{U}(t,k))\leq0
\end{equation}
holds for  $(t,k)\in(0,\infty)\times\mathbb{R}^3.$
\end{theorem}
\noindent \emph{Proof.} We will consider the linearized homogeneous
system \eqref{4.1} in Fourier space. For this purpose, by taking
Fourier transform in $x$ for the linearized homogeneous system
\eqref{4.1}, then $\hat{U}=[\hat{\rho},~ \hat{u},~ \hat{\Theta},~
\hat{E},~ \hat{B}]$ satisfies
\begin{equation}
\label{4.6}
\left\{\begin{aligned}
&\partial_t \hat{\rho} + ik\cdot \hat{u}=0, \\
&\partial_t \hat{u} + ik \hat{\rho} +ik \hat{\Theta}+\hat{E}+ \hat{u}   =0, \\
&\partial_t \hat{\Theta}+\frac{2}{3}ik\cdot \hat{u} + \hat{\Theta}= 0 , \\
&  \partial_t \hat{E}-ik\times \hat{B}-\hat{u} =0, \\
&  \partial_t \hat{B}+ik\times \hat{E}=0,\\
& ik\cdot \hat{E}=-\hat{\rho} ,\quad ik\cdot \hat{B}=0,\quad
(t,k)\in(0,\infty)\times\mathbb{R}^3,
\end{aligned} \right.
\end{equation}
Firstly, one can acquire from the first five equations of the system
\eqref{4.6} that
\begin{equation}\label{4.7}
\partial_t\left|\left[\hat{\rho},~\hat{u},~{\frac{\sqrt 6}{2}}
\hat{\Theta},~\hat{E},~\hat{B}\right]\right|^2
+2|\hat{u}|^2+3|\hat{\Theta}|^2=0.
\end{equation}
Multiplying the second equation of the system \eqref{4.6} by
$\overline{ik\hat{\rho}} $, utilizing integration by parts in $t$
and replacing $\partial_t\hat{\rho}$ by the first equation of the
system \eqref{4.6}, we have
\begin{equation}\label{A1}
\partial_t(\hat{u} ~|ik\hat{\rho} )+(1+|k|^2)|\hat{\rho}|^2=|k\cdot\hat{u}|^2
-|k|^2\hat{\Theta}~\overline{\hat{\rho}}+ik\cdot\hat{u}~\overline{\hat{\rho}}
.\end{equation}
 Multiplying the second
equation of the system \eqref{4.6} by $\overline{ik\hat{\Theta}} $,
utilizing integration by parts in $t$ and replacing
$\partial_t\hat{\Theta}$ by the third equation of the system
\eqref{4.6}, we also have
\begin{equation}\notag
\partial_t(\hat{u} ~|ik\hat{\Theta} )+|k|^2 |\hat{\Theta}|^2
=\frac{2}{3}|k\cdot\hat{u}|^2
+2ik\cdot\hat{u}~\overline{\hat{\Theta}}-|k|^2\hat{\rho}~\overline{\hat{\Theta}}
-\hat{\rho}~\overline{\hat{\Theta}},
\end{equation}
putting it together with \eqref{A1} and taking the real part after
utilizing the Cauchy-Schwarz inequality, one has
 $$\partial_t  \mathcal {R}\left(\hat{u}~|ik\left(\hat{\rho}+\hat{\Theta}\right) \right)
 +\gamma |\hat{\rho}|^2
 \leq C |k\cdot\hat{u} |^2+C|\hat{\Theta}|^2.$$
Multiplying it by $\displaystyle\frac{1}{1+|k|^2} $, one can obtain
\begin{equation}\label{4.8}
\frac{\partial_t  \mathcal
{R}\left(\hat{u}~|ik\left(\hat{\rho}+\hat{\Theta}\right)
\right)}{1+|k|^2}
 +\frac{\gamma |\hat{\rho}|^2}{1+|k|^2}
 \leq C(|\hat{u} |^2+|\hat{\Theta}|^2).
 \end{equation}
Similarly,  multiplying the second equation of the system
\eqref{4.6} by $\overline{\hat{E}}$, utilizing integration by parts
in $t$ and replacing $\partial_t \hat{E}$ by the fourth equation of
the system \eqref{4.6}, we have
\begin{equation} \label{4.9}
 {\partial_t}\left( {\hat u\left| {\hat E} \right.} \right) + \left( {{{\left|
{\hat E} \right|}^2} + {{\left| {k \cdot \hat E} \right|}^2}}
\right) = {\left| {\hat u} \right|^2} - \hat \Theta ~\overline{ \hat
\rho} - ik \times \overline{ \hat B }\cdot \hat u - \hat u \cdot
 \overline{\hat E}.\end{equation}
Using the Cauchy-Schwarz inequality and multiplying it by
$\displaystyle{\frac{|k|^2}{\left(1+|k|^2\right)^2}}$, one can
obtain
\begin{equation} \label{4.10}
\begin{split}
 {\partial_t}\frac{|k|^2\mathcal {R}\left( {\hat u\left| {\hat E} \right.} \right)}{\left(1+|k|^2\right)^2}
 & + \frac{\gamma|k|^2\left( {{{\left|
{\hat E} \right|}^2}  + {{\left| {k \cdot \hat E} \right|}^2}}
\right)}{\left(1+|k|^2\right)^2}\\
& \leq C \left(   {\frac{\left| {\hat \rho} \right|^2}{1+|k|^2}}
 +\left|
{\hat u} \right|^2 +\left| {\hat\Theta} \right|^2
 \right)+\frac{|k|^2 \mathcal {R}\left(- ik \times \overline{ \hat B }\cdot
\hat u \right)}{\left(1+|k|^2\right)^2}.
\end{split}
\end{equation}
Similarly, from the fourth and fifth equations of the system
\eqref{4.6}, one has
$$
{\partial _t}\left( { - ik \times \hat B\left| {\hat E} \right.}
\right) + {\left| {k \times \hat B} \right|^2} = {\left| {k \times
\hat E} \right|^2} - \left( {ik \times \hat B\left| {\hat u}
\right.} \right),
$$
which after utilizing Cauchy-Schwarz inequality and multiplying it
by $\displaystyle\frac{1}{(1+|k|^2)^2}$ , gives
\begin{equation}\label{4.11}
\partial_t\frac{\mathcal {R}(-ik\times\hat{B}|\hat{E})}{(1+|k|^2)^2}+\gamma
\frac{|k\times\hat{B}|^2}{(1+|k|^2)^2}\leq
\frac{|k|^2|\hat{E}|^2}{(1+|k|^2)^2}+C|\hat{u}|^2.
\end{equation}
Lastly, we define the time-frequency Lyapunov functional as
\begin{equation}\notag
\begin{split}
\mathcal
{E}(\hat{U}(t,k))=\left|\left[\hat{\rho},~\hat{u},~{\frac{\sqrt
6}{2}} \hat{\Theta},~\hat{E},~\hat{B}\right]\right|^2 &+\mathcal
{K}_1\frac{  \mathcal
{R}\left(\hat{u}~|ik\left(\hat{\rho}+\hat{\Theta}\right)
\right)}{1+|k|^2}\\
&+\mathcal {K}_2\frac{|k|^2\mathcal {R}\left( {\hat u\left| {\hat E}
\right.} \right)}{\left(1+|k|^2\right)^2} +\mathcal
{K}_3\frac{\mathcal {R}(-ik\times\hat{B}|\hat{E})}{(1+|k|^2)^2},
\end{split}
\end{equation}
where, constants  $0<\mathcal {K}_3\ll\mathcal {K}_2\ll\mathcal
{K}_1\ll 1$ are to be determined later. Notice that as soon as
$0<\mathcal {K}_j \ll1 _{(j=1,2,3)}$ be sufficiently small, then
\eqref{4.4} holds true. Moreover, by setting $0<\mathcal
{K}_3\ll\mathcal {K}_2\ll\mathcal {K}_1\ll 1$ be sufficiently small
with $\mathcal {K}_2^{\frac{3}{2}}\ll\mathcal {K}_3$, taking the
summation of \eqref{4.7}, \eqref{4.8}$\times\mathcal {K}_1$,
\eqref{4.10}$\times\mathcal {K}_2$ and \eqref{4.11}$\times\mathcal
{K}_3$, one has
\begin{equation}\label{4.12}
\partial_t\mathcal {E}(\hat{U}(t,k))+ \frac{\gamma\left| {\hat \rho} \right|^2}{1+|k|^2}
+\gamma\left(|\hat{u}|^2+|\hat{\Theta}|^2\right)+\frac{\gamma
|k|^2}{(1+|k|^2)^2}|[ \hat{E},\hat{B}]|^2\leq0,
\end{equation}
where we have used the Cauchy-Schwarz inequality as follows
$$\mathcal
{K}_2\frac{|k|^2\mathcal
{R}(-ik\times\overline{\hat{B}}\cdot\hat{u})}{(1+|k|^2)^2} \leq
\frac{\mathcal {K}_2^{\frac{1}{2}}|k|^4|\hat{u}|^2}{2(1+|k|^2)^2}+
\frac{\mathcal
{K}_2^{\frac{3}{2}}|k|^2||\hat{B}|^2}{2(1+|k|^2)^2}.$$ Then, by
noticing $\mathcal {E}(\hat{U})\sim |\hat{U}|^2$ and
$$\frac{\gamma |k|^2}{(1+|k|^2)^2}|\hat{U}|^2\leq
\frac{\gamma\left| {\hat \rho} \right|^2}{1+|k|^2}
+\gamma\left(|\hat{u}|^2+|\hat{\Theta}|^2\right)+\frac{\gamma
|k|^2}{(1+|k|^2)^2}|[ \hat{E},\hat{B}]|^2,$$ one can obtain
\eqref{4.5} from \eqref{4.12}. Now, we have finished the proof of Theorem \ref{thm4.1}.\hfill$\Box$\\
\\
\indent It is straightforward from Theorem \ref{thm4.1} to obtain
the pointwise time-frequency estimate on the norm of
$|\widehat{U}(t,k)|$ in terms of the given initial data norm of
$|\widehat{U}_0(k)|$.
\begin{coro}\label{Corollary4.1}
Suppose $U(t,x),(t,x)\in(0,\infty)\times\mathbb{R}^3 $ to be a
 solution to the initial problem \eqref{4.1}-\eqref{4.3}.
Then, there exist constants $\gamma>0,C>0$ such that
\begin{equation}\label{4.13}
|\widehat{U}(t,k)|\leq Ce^{-\frac{\gamma
|k|^2t}{(1+|k|^2)^2}}|\widehat{U}_0(k)|
\end{equation}
holds for $t\geq0$ and $k\in{\mathbb R}^3$.
\end{coro}
\vspace{0.3cm}

\noindent 4.2 \textbf{$L^p-L^q$ time-decay property.}

Formally, the
solution to the initial problem \eqref{4.1}-\eqref{4.2} is presented
as
\begin{equation}\label{4.14}
U(t)=e^{tL}U_0,
\end{equation}
here, $e^{tL}$ is named as the linear solution operator for
${t\geq0}$. The main result of this subsection, whose proof will be
omitted here, is stated as follows; see
\cite{Duan11},\cite{DH95},\cite{SK84}.
\begin{theorem}\label{thm4.2}
Let $1\leq p,r\leq 2\leq q\leq\infty,l\geq0$ and an integer
$m\geq0$. Define
\begin{equation}\label{4.15}
{[l + 3(\frac{1}{r} - \frac{1}{q})]_ + } = \left\{
{\begin{array}{*{20}{c}}
   l & { \mbox{if} ~~  r = q = 2\ \mbox{ and }{\rm{ }}l{\mbox{ is an integer,}}}  \\
   {{{[l + 3(\frac{1}{r} - \frac{1}{q})]}_ - } + 1} & \mbox{otherwise,}  \\
\end{array}} \right.
\end{equation}
where, we use $[\cdot]_-$ to denote the integer part of the
argument. Assume $U_0$ satisfies \eqref{4.3}. Then, $e^{tL}$
satisfies the following time-decay property:
\begin{equation}\label{4.16}
\begin{split}
 \left\|{\nabla^me^{tL}U_0}\right\|_{L^q}
\leq
C(1+t)^{-\frac{3}{2}(\frac{1}{p}-\frac{1}{q})-\frac{m}{2}}\left\|{U_0}\right\|_{L^p}
+C(1+t)^{-\frac{t}{2}}\left\|{\nabla^{m+[l + 3(\frac{1}{r} -
\frac{1}{q})]_ +}U_0}\right\|_{L^r}
\end{split}
\end{equation}
for $t\geq0$, here $C=C(p,q,r,l,m)$.
\end{theorem}
\vspace{0.3cm}

  \noindent 4.3 \textbf{ Representation of solutions.} We investigate the explicit solution
$U=[\rho,u,\Theta,E,B]=e^{tL}U_0$ to the initial problem
\eqref{4.1}-\eqref{4.2} with the compatible condition \eqref{4.3} or
equivalently the equations \eqref{4.6} in Fourier space in this
subsection. We find that $\rho$, $\Theta$, $\nabla\cdot u$ satisfy
the same equation which is of order three and is different from that
of the isentropic Euler-Maxwell system. The main purpose is to prove
Theorem \ref{thm4.3} presented at the bottom of this subsection.

From the first three equations of \eqref{4.1} and $\nabla\cdot
E=-\rho,$ one has
\begin{equation}\label{4.18}
\partial_{ttt}\rho
+2\partial_{tt}\rho-\frac{5}{3}\partial_{t}\triangle\rho+2\partial_{t}\rho+\rho-\triangle\rho=0,
\end{equation}
with initial data
\begin{equation}\label{4.19}
\left\{\begin{split}
&\rho|_{t=0}=\rho_0=-\nabla\cdot E_0,\\
&\partial_{t}\rho|_{t=0}=-\nabla\cdot u_0,\\
&\partial_{tt}\rho|_{t=0}=\nabla\cdot
u_0+\triangle\rho_0-\rho_0+\triangle\Theta_0.
\end{split}
\right.\end{equation} After taking Fourier transform in $x$ for
\eqref{4.18} and \eqref{4.19}, one has
\begin{equation}\label{4.20}
\partial_{ttt}\hat{\rho}
+2\partial_{tt}\hat{\rho}+\left(2+\frac{5}{3}|k|^2\right)\partial_{t}\hat{\rho}
+\left(1+|k|^2\right)\hat{\rho}=0,
\end{equation}
with initial data
\begin{equation}\label{4.21}
\left\{\begin{split}
&\hat{\rho}|_{t=0}=\hat{\rho}_0=-ik\cdot \hat{E}_0,\\
&\partial_{t}\hat{\rho}|_{t=0}=-ik\cdot \hat{u}_0,\\
&\partial_{tt}\hat{\rho}|_{t=0}=-\left(1+|k|^2\right)\hat{\rho}_0+ik\cdot
\hat{u}_0-|k|^2\Theta_0.
\end{split}
\right.\end{equation} The characteristic equation of \eqref{4.20} is
$$ F(\mathcal {X}):=\mathcal {X}^3+2\mathcal {X}^2+\left(2+\frac{5}{3}|k|^2\right)\mathcal {X}+1+|k|^2=0.$$
For the roots of the previous characteristic equation and their
 properties, we obtain
\begin{lemma}\label{L4.1}
Suppose $|k|\neq0.$ Then, $F(\mathcal {X})=0,$ $\mathcal
{X}\in\mathbb{C}$ has a real root
$\sigma=\sigma(|k|)\in(-1,-\frac{3}{5})$ and two conjugate complex
roots $\mathcal {X}_\pm=\beta\pm i\omega$ with
$\beta=\beta(|k|)\in(-\frac{7}{10},-\frac{1}{2})$ and
$\omega=\omega(|k|)\in({\frac{\sqrt6}{3}},+\infty)$ which satisfy
\begin{equation}\label{4.22}
\beta=-1-\frac{\sigma}{2},~
\omega=\frac{1}{2}\sqrt{3\sigma^2+4\sigma+4+\frac{20}{3}|k|^2}.
\end{equation}
$\sigma,\beta,\omega$ are smooth in $|k|>0$, and $\sigma(|k|)$ is
strictly increasing over $|k|>0$, with
$$ \lim_{|k|\longrightarrow0}\sigma(|k|)=-1,~\lim_{|k|\longrightarrow\infty}\sigma(|k|)=-\frac{3}{5}.$$
Furthermore, the asymptotic behavior in the following hold true:
$$\sigma(|k|)=-1+O(1)|k|^2,~\beta(|k|)=-\frac{1}{2}-O(1)|k|^2,~\omega(|k|)=\frac{\sqrt{3}}{2}+O(1)|k| $$
whenever $|k|\leq1$ is sufficiently small, and
$$\sigma(|k|)=-\frac{3}{5}-O(1)|k|^{-2},~\beta(|k|)=-\frac{7}{10}+O(1)|k|^{-2},~\omega(|k|)=O(1)|k|$$
whenever $|k|\geq1$ is sufficiently large. Here and in the sequel,
we use $O(1)$ to denote a generic strictly positive constant.
\end{lemma}
\noindent \emph{Proof.} Let $|k|\neq0.$ Firstly, we search the
possibly existing real root for equation $F(\mathcal {X})=0$ in
$\mathcal {X}\in \mathbb{R}$. From that
 $$F'(\mathcal {X})=3\mathcal {X}^2+4\mathcal {X}+2+\frac{5}{3}|k|^{2}=
 3(\mathcal {X}+\frac{2}{3})^2+\frac{2}{3}+\frac{5}{3}|k|^{2}>0,$$
 and
 $F(-1)=-\frac{2}{3}|k|^{2}<0,~F(-\frac{3}{5})=\frac{38}{125}>0,$
 one can obtain that equation $F(\mathcal {X})=0$ indeed has one and only one real
 root denoted by $\sigma=\sigma(|k|)$ which satisfies $-1<\sigma<-\frac{3}{5}.$
  After taking derivative of $F(\sigma(|k|))=0$ in $|k|$, one has
 $$\sigma'(|k|)=\frac{-|k|\left( 2+\frac{10}{3}\sigma\right)}{3\sigma^2+4\sigma+2+\frac{5}{3}|k|^{2}}>0,$$
 so that $\sigma(\cdot)$ is strictly increasing over $|k|>0.$ Since
 $F(\sigma)=0$ can be represented as
 $$\sigma\left[ \frac{\sigma(\sigma+2)}{2+\frac{5}{3}|k|^{2}}+1\right]=-\frac{1+|k|^{2}}{2+\frac{5}{3}|k|^{2}},$$
then $\sigma$ has limits $-1$ and $-\frac{3}{5}$ as $|k|\rightarrow
0$ and $|k|\rightarrow \infty$, respectively.

$F(\sigma(|k|))=0$ is also equivalent with
$$\sigma+1=\frac{\frac{2}{3}|k|^{2}+(\sigma+1)^2}{(\sigma+1)^2+1+\frac{5}{3}|k|^{2}}$$
or
$$\sigma+\frac{3}{5}=\frac{-\frac{1}{5}+\frac{3}{5}\sigma(\sigma+2)}{\sigma(\sigma+2)+2+\frac{5}{3}|k|^{2}}.$$
Therefore, it follows that $\sigma(|k|)=-1+O(1)|k|^2$ whenever
$|k|<1$ is sufficiently small and
$\sigma(|k|)=-\frac{3}{5}-O(1)|k|^{-2} $ whenever $|k|\geq 1$ is
sufficiently large. Next, let us search roots of $F(\mathcal {X})=0$
in $\mathcal {X}\in \mathbb{C}$. Since $F(\sigma)=0$ with $\sigma\in
\mathbb{R}$, $F(\mathcal {X})=0$ can be decomposed as
$$ F(\mathcal {X})=(\mathcal {X}-\sigma)\left[\left(\mathcal {X}+1+\frac{\sigma}{2} \right)^2+\frac{3}{4}\sigma^2
+\sigma+\frac{5}{3}|k|^2+1 \right]=0.$$ Then, there exist two
conjugate complex roots $\mathcal {X}_\pm=\beta\pm i\omega$ which
satisfy
$$ \left(\mathcal {X}+1+\frac{\sigma}{2} \right)^2+\frac{3}{4}\sigma^2+\sigma+\frac{5}{3}|k|^2+1=0.$$
By solving the previous equation, one can obtain that
$\beta=\beta(|k|),~\omega=\omega(|k|)$ take the form of
\eqref{4.22}. It is straightforward from the asymptotic behavior of
$\sigma(|k|)$ at $|k|=0$ and $\infty$ to obtain that of
$\beta(|k|),~\omega(|k|)$. Then, we have finished the proof of Lemma
\ref{L4.1}.\hfill$\Box$

Based on Lemma \ref{L4.1}, we define the solution of \eqref{4.20} as
\begin{equation}\label{4.23}
\hat{\rho}(t,k)=c_1(k)e^{\sigma t}+e^{\beta t}\left(
c_2(k)\cos\omega t+c_3(k)\sin\omega t \right),
\end{equation}
here $c_j(k),~1\leq j \leq 3,$ is to be chosen by \eqref{4.21}
later. Again using $\nabla\cdot E=-\rho,$ \eqref{4.23} implies
\begin{equation}\label{4.24}
\tilde{k}\cdot \hat{E}(t,k)=i|k|^{-1}\left( c_1(k)e^{\sigma
t}+e^{\beta t}\left( c_2(k)\cos\omega t+c_3(k)\sin\omega t
\right)\right).
\end{equation}
Here and in the sequel $\tilde{k}=\frac{k}{|k|}.$ In fact,
\eqref{4.23} implies
\begin{equation}\label{4.25}
\left[ {\begin{array}{*{20}{c}}
   {\hat \rho {|_{t = 0}}}  \\
   {{\partial _t}\hat \rho {|_{t = 0}}}  \\
   {{\partial _{tt}}\hat \rho {|_{t = 0}}}  \\

 \end{array} } \right] = A\left[ {\begin{array}{*{20}{c}}
   {{c_1}}  \\
   {{c_2}}  \\
   {{c_3}}  \\

 \end{array} } \right],\quad A = \left[ {\begin{array}{*{20}{c}}
   1 & 1 & 0  \\
   \sigma  & \beta  & \omega   \\
   {{\sigma ^2}} & {{\beta ^2} - {\omega ^2}} & {2\beta \omega }  \\

 \end{array} } \right].\end{equation}
It is directly to obtain that
$$\det A = \omega \left[ {{\omega ^2} + {{\left( {\sigma  - \beta } \right)}^2}} \right]
= \omega \left( {3{\sigma ^2} + 4\sigma  + 2 + \frac{5} {3}{{\left|
k \right|}^2}} \right) > 0
$$
and
\[{A^{ - 1}} = \frac{1}
{{\det A}}\left[ {\begin{array}{*{20}{c}}
   {\left( {{\beta ^2} + {\omega ^2}} \right)\omega } & { - 2\beta \omega } & \omega   \\
   {\sigma \left( {\sigma  - 2\beta } \right)\omega } & {2\beta \omega } & { - \omega }  \\
   {\sigma \left( {{\beta ^2} - {\omega ^2} - \sigma \beta } \right)} & {{\omega ^2} + {\sigma ^2}
   - {\beta ^2}} & {\beta  - \sigma }  \\

 \end{array} } \right].\]
From \eqref{4.25} and \eqref{4.21}, one has
\begin{equation}\notag
\begin{split}
[{c_1},~ &   {c_2},~   {c_3}]^T =\frac{1} {{3{\sigma ^2} +
4\sigma  + 2 + \frac{5} {3} {{\left| k \right|}^2}}}\\
&\left[ {\begin{array}{*{20}{c}}
   {{\beta ^2} + {\omega ^2} - \left( {1 + {{\left| k \right|}^2}} \right)} & { - i\left| k \right|\left( {2\beta  + 1} \right)} & { - {{\left| k \right|}^2}}  \\
   {{\sigma ^2} - 2\sigma \beta  + \left( {1 + {{\left| k \right|}^2}} \right)} & {i\left| k \right|\left( {2\beta  + 1} \right)} & {{{\left| k \right|}^2}}  \\
   {\frac{{\sigma \left( {{\beta ^2} - {\omega ^2} - \sigma \beta } \right) - \left( {\beta  - \sigma } \right)\left( {1 + {{\left| k \right|}^2}} \right)}}
{\omega }} & {\frac{{i\left| k \right|}} {\omega }\left( {{\beta ^2}
- {\sigma ^2} - {\omega ^2} + \beta  - \sigma } \right)} &
{\frac{{\sigma  - \beta }}
{\omega }{{\left| k \right|}^2}}  \\
 \end{array} } \right]\left[ {\begin{array}{*{20}{c}}
   {{{\hat \rho }_0}}  \\
   {\tilde k \cdot {{\hat u}_0}}  \\
   {{{\hat \Theta }_0}}  \\
 \end{array} } \right].
 \end{split}\end{equation}
Where, we use $[\cdot]^T$  to denote the transpose of a vector.
Making further simplifications with the form of $\beta$ and
$\omega$, we have
\begin{equation}\label{4.26}
\begin{split}
[&{c_1},~   {c_2},~   {c_3}]^T =\frac{1} {{3{\sigma ^2} +
4\sigma  + 2 + \frac{5} {3} {{\left| k \right|}^2}}}\\
&\left[ {\begin{array}{*{20}{c}}
   {{{\left( {\sigma  + 1} \right)}^2} + \frac{2}
{3}{{\left| k \right|}^2}} & { - i\left| k \right|\left( {\sigma  + 1} \right)} & { - {{\left| k \right|}^2}}  \\
   {2{\sigma ^2} + \sigma  + 1 + {{\left| k \right|}^2}} & {i\left| k \right|\left( {\sigma  + 1} \right)} & {{{\left| k \right|}^2}}  \\
   {\frac{{{\sigma ^2} + \frac{3}
{2}\sigma  + \left( {1 + {{\left| k \right|}^2}} \right) - \frac{1}
{6}\sigma {{\left| k \right|}^2}}} {\omega }} & {\frac{{i\left| k
\right|}} {\omega }\left( {\frac{3} {2}{\sigma ^2} + \frac{3}
{2}\sigma  + 1 + \frac{5} {3}{{\left| k \right|}^2}} \right)} &
{\frac{{1 + \frac{3} {2}\sigma }}
{\omega }{{\left| k \right|}^2}}  \\
 \end{array} } \right]\left[ {\begin{array}{*{20}{c}}
   {{{\hat \rho }_0}}  \\
   {\tilde k \cdot {{\hat u}_0}}  \\
   {{{\hat \Theta }_0}}  \\
 \end{array} } \right].
 \end{split}\end{equation}
Similarly, from the first three equations of \eqref{4.1} and
$\nabla\cdot E=-\rho,$ one has
\begin{equation}\label{4.27}
\partial_{ttt}\hat{\Theta}
+2\partial_{tt}\hat{\Theta}+\left(2+\frac{5}{3}|k|^2\right)\partial_{t}\hat{\Theta}
+\left(1+|k|^2\right)\hat{\Theta}=0,
\end{equation}
with initial data
\begin{equation}\label{4.28}
\left\{\begin{split}
&\hat{\Theta}|_{t=0}=\hat{\Theta}_0,\\
&\partial_{t}\hat{\Theta}|_{t=0}=\frac{2}{3}ik\cdot \hat{u}_0,\\
&\partial_{tt}\hat{\Theta}|_{t=0}=-\frac{2}{3}\left(1+|k|^2\right)\hat{\rho}_0+\frac{4}{3}ik\cdot
\hat{u}_0+\left(1-\frac{2}{3}|k|^2\right)\hat{\Theta}_0.
\end{split}
\right.\end{equation} From Lemma \ref{L4.1}, one can define the
solution of \eqref{4.27} as
\begin{equation}\label{4.29}
\hat{\Theta}(t,k)=c_4(k)e^{\sigma t}+e^{\beta t}\left(
c_5(k)\cos\omega t+c_6(k)\sin\omega t \right),
\end{equation}
where $c_j(k),~4\leq j \leq 6,$ is to be chosen by \eqref{4.28}
later. In fact, \eqref{4.29} implies
\begin{equation}\label{4.30}
\begin{split}
[&{c_4},~   {c_5},~   {c_6}]^T =\frac{1} {{3{\sigma ^2} +
4\sigma  + 2 + \frac{5} {3} {{\left| k \right|}^2}}} \\
& \left[ {\begin{array}{*{20}{c}}
   {{\sigma ^2} + 2\sigma  + \frac{4}
{3} + \frac{2} {3}{{\left| k \right|}^2}} & {\frac{4} {3}\left| k
\right|\left( {2 + \frac{\sigma } {2}} \right)i} & {1 - \frac{2}
{3}{{\left| k \right|}^2}}  \\
   {2{\sigma ^2} + 3\sigma  + \frac{8}
{3} + \frac{2} {3}{{\left| k \right|}^2}} & { - \frac{4} {3}\left| k
\right|\left( {2 + \frac{\sigma } {2}} \right)i} & { - \left( {1 -
\frac{2}
{3}{{\left| k \right|}^2}} \right)}  \\
   {\frac{{ - \frac{1}
{2}{\sigma ^2} + \sigma  - \frac{2} {3}\sigma {{\left| k \right|}^2}
+ \frac{2} {3} - {{\left| k \right|}^2}}} {\omega }} &
{\frac{{\frac{3} {2}{\sigma ^2} - 3\sigma  - 2 + \frac{5}
{3}{{\left| k \right|}^2}}} {\omega }} & {\frac{{ - 1 - \frac{3}
{2}\sigma  + \frac{2} {3}{{\left| k \right|}^2} + \sigma {{\left| k
\right|}^2}}}
{\omega }}  \\
 \end{array} } \right]\left[ {\begin{array}{*{20}{c}}
   {{{\hat \rho }_0}}  \\
   {\tilde k \cdot {{\hat u}_0}}  \\
   {{{\hat \Theta }_0}}  \\
 \end{array} } \right].
 \end{split}\end{equation}
Similarly, again from the first three equations of \eqref{4.1} and
$\nabla\cdot E=-\rho,$ one also has
\begin{equation}\label{4.31}
\partial_{ttt}{(\tilde{k}\cdot \hat{u})}
+2\partial_{tt}{(\tilde{k}\cdot
\hat{u})}+\left(2+\frac{5}{3}|k|^2\right)\partial_{t}{(\tilde{k}\cdot
\hat{u})} +\left(1+|k|^2\right){(\tilde{k}\cdot \hat{u})}=0.
\end{equation}
Initial data is given by
\begin{equation}\label{4.32}
\left\{\begin{split}
&{(\tilde{k}\cdot \hat{u})}|_{t=0}={(\tilde{k}\cdot \hat{u})}_0,\\
&\partial_{t}{(\tilde{k}\cdot \hat{u})}|_{t=0}=\frac{{ - i}}
{{\left| k \right|}}\left( {1 + {{\left| k \right|}^2}} \right){\hat
\rho _0} - \tilde k \cdot {\hat u_0} - i\left| k \right|{\hat \Theta
_0},
\\
&\partial_{tt}{(\tilde{k}\cdot \hat{u})}|_{t=0}=\frac{i} {{\left| k
\right|}}\left( {1 + {{\left| k \right|}^2}} \right){\hat \rho _0} -
\frac{5} {3}{\left| k \right|^2}\tilde k \cdot {\hat u_0} + 2i\left|
k \right|{\hat \Theta _0}.
\end{split}
\right.\end{equation} After tenuous computation, one has
\begin{equation}\label{4.33}
{\tilde{k}\cdot \hat{u}}(t,k)=c_7(k)e^{\sigma t}+e^{\beta  t}\left(
c_8(k)\cos\omega t+c_9(k)\sin\omega t \right),
\end{equation}
with
\begin{equation}\label{4.34}
\begin{split}
[{c_7},~   & {c_8},~   {c_9}]^T =\frac{1} {{3{\sigma ^2} +
4\sigma  + 2 + \frac{5} {3} {{\left| k \right|}^2}}} \\
&   \left(~ {\left[ {\begin{array}{*{20}{c}}
   {{\sigma ^2} + 2\sigma  + 2 + \frac{5}
{3}{{\left| k \right|}^2}} & { - 2 - \sigma  - \frac{5}
{3}{{\left| k \right|}^2}} & { - \sigma \left| k \right|i}  \\
   {2{\sigma ^2} + 2\sigma } & {2 + \sigma  + \frac{5}
{3}{{\left| k \right|}^2}} & {\sigma \left| k \right|i}  \\
   {\frac{{{\sigma ^2} - \frac{5}
{3}\sigma {{\left| k \right|}^2}}} {\omega }} & {\frac{{\frac{3}
{2}{\sigma ^2} + \frac{5} {2}\sigma {{\left| k \right|}^2}}} {\omega
}} & {\frac{{ - \frac{3} {2}{\sigma ^2} - 3\sigma  - 2 - \frac{5}
{3}{{\left| k \right|}^2}}}
{\omega }}  \\
 \end{array} } \right]} \right. \\
& + \left. {\left[ {\begin{array}{*{20}{c}}
   { - \left( {1 + \sigma } \right){{\left| k \right|}^{ - 1}}\left( {1 + {{\left| k \right|}^2}} \right)i} & 0 & 0  \\
   {\left( {1 + \sigma } \right){{\left| k \right|}^{ - 1}}\left( {1 + {{\left| k \right|}^2}} \right)i} & 0 & 0  \\
   { - \left( {\frac{3}
{2}{\sigma ^2} + \frac{3} {2}\sigma  + 1 + \frac{5}
{3}{{\left| k \right|}^2}} \right){{\left| k \right|}^{ - 1}}\left( {1 + {{\left| k \right|}^2}} \right)i} & 0 & 0  \\
 \end{array} } \right]} ~\right)   \left[ {\begin{array}{*{20}{c}}
   {{{\hat \rho }_0}}  \\
   {\tilde k \cdot {{\hat u}_0}}  \\
   {{{\hat \Theta }_0}}  \\
 \end{array} } \right].
 \end{split}\end{equation}

Next, for $(t,k)\in(0,\infty)\times\mathbb{R}^3 $, let us solve
\begin{equation}\notag
\left\{
\begin{split}
& M_{1}(t,k):=-\tilde{k}\times(\tilde{k}\times \hat{u}(t,k)),\\
& M_{2}(t,k):=-\tilde{k}\times(\tilde{k}\times \hat{E}(t,k)),\\
& M_{3}(t,k):=-\tilde{k}\times(\tilde{k}\times \hat{B}(t,k)).
\end{split}
\right.
\end{equation}
Taking the curl for the second, fourth and fifth equations of the
system \eqref{4.1}, one has
\begin{equation}\notag
\left\{
\begin{split}
& \partial_t(\nabla\times u)+\nabla\times E+\nabla\times u=0,\\
& \partial_t(\nabla\times E)-\nabla\times(\nabla\times B)-\nabla\times u=0,\\
& \partial_t(\nabla\times B)+\nabla\times(\nabla\times E)=0.\\
\end{split}
\right.
\end{equation}
Taking Fourier transform in $x$ for the previous system, it follows
that
\begin{equation}\notag
\left\{
\begin{split}
& \partial_tM_{1}=-M_{1}-M_2,\\
& \partial_tM_2=M_{1}+ik\times M_3,\\
& \partial_tM_3=-ik\times M_2.\\
\end{split}
\right.
\end{equation}
Initial data is given as
\begin{equation}\notag
\begin{split}
[M_{1}, M_{2}, M_{3}]|_{t=0}=[M_{1,0}, M_{2,0}, M_{3,0}].
\end{split}
\end{equation}
Here,
\begin{equation}\notag
M_{1,0}=-\tilde{k}\times(\tilde{k}\times \hat{u}_{0}),~
M_{2,0}=-\tilde{k}\times(\tilde{k}\times
\hat{E}_{0}),~M_{3,0}=-\tilde{k}\times(\tilde{k}\times \hat{B}_{0}).
\end{equation}
It straightforward to get
\begin{equation}\label{4.35}
\partial_{ttt}M_2
+\partial_{tt}M_2+\left(1+|k|^2\right)\partial_{t}M_2 +|k|^2M_2=0,
\end{equation}
with initial data
\begin{equation}\label{4.36}
\left\{
\begin{split}
& M_2|_{t=0}=M_{2,0},\\
& \partial_tM_2|_{t=0}=M_{1,0}+ik\times M_{3,0},\\
& \partial_{tt}M_2|_{t=0}=-M_{1,0}-(1+|k|^2)M_{2,0}.\\
\end{split}
\right.
\end{equation}
The characteristic equation of \eqref{4.35} is
$$ F_*(\mathcal {X}):=\mathcal {X}^3+\mathcal {X}^2+\left(1+|k|^2\right)\mathcal {X}+|k|^2=0.$$
For the roots of the previous characteristic equation and their
 properties, we have
\begin{lemma}\label{L4.2}
Suppose $|k|\neq0.$ Then, $F_*(\mathcal {X})=0,$ $\mathcal
{X}\in\mathbb{C}$ has a real root $\sigma_*=\sigma_*(|k|)\in(-1,0)$
and two conjugate complex roots $\mathcal {X}_\pm=\beta_*\pm
i\omega_*$ with $\beta_*=\beta_*(|k|)\in(-\frac{1}{2},0)$ and
$\omega_*=\omega_*(|k|)\in({\frac{\sqrt 6}{3}},+\infty)$ which
satisfy
\begin{equation}\label{4.37}
\beta_*=-\frac{1}{2}-\frac{\sigma_*}{2},~
\omega_*=\frac{1}{2}\sqrt{3\sigma_*^2+2\sigma_*+3+4|k|^2}.
\end{equation}
$\sigma_*,\beta_*,\omega_*$ are smooth in $|k|>0$, and
$\sigma_*(|k|)$ is strictly decreasing over $|k|>0$, with
$$ \lim_{|k|\longrightarrow0}\sigma_*(|k|)=0,~\lim_{|k|\longrightarrow\infty}\sigma_*(|k|)=-1.$$
Furthermore, the asymptotic behavior as follows hold true:
$$\sigma_*(|k|)=-O(1)|k|^2,~\beta_*(|k|)=-\frac{1}{2}+O(1)|k|^2,~\omega_*(|k|)=\frac{\sqrt{3}}{2}+O(1)|k| $$
whenever $|k|\leq1$ is sufficiently small, and
$$\sigma_*(|k|)=-1+O(1)|k|^{-2},~\beta_*(|k|)=-O(1)|k|^{-2},~\omega_*(|k|)=O(1)|k|$$
whenever $|k|\geq1$ is sufficiently large.
\end{lemma}
Similarly as before, we obtain
\begin{equation}\label{4.38}
\begin{split}
 M_{1} (t,k)
  &=  - \frac{{c_{10} (k)}}{{1 + \sigma_* }}e^{\sigma_* t}  - \frac{{c_{11} (k)}}
 {{(1 + \beta_* )^2  + \omega_* ^2 }}e^{\beta_* t} \big[(1 + \beta_* )\cos \omega_* t
+ \omega_* \sin \omega_* t\big]\\
&\quad - \frac{{c_{12} (k)}}{{(1 + \beta_* )^2  + \omega_* ^2
}}e^{\beta_* t}
 \left[(1 + \beta_* )\sin \omega_* t - \omega_* \cos \omega_* t\right],
 \end{split}
\end{equation}
\begin{equation}\label{4.39}
M_{2}(t,k)=c_{10}(k)e^{\sigma_*  t}+e^{\beta_*
t}\left[c_{11}(k)\cos\omega_*  t+c_{12}(k)\sin \omega_*  t\right],
\end{equation}
and
\begin{equation}\label{4.40}
\begin{split}
 M_3 (t,k)
  &= - ik \times \frac{{c_{10} (k)}}{\sigma_*  }e^{\sigma_*  t}  - ik \times \frac{{c_{11} (k)}}
 {{\beta_*  ^2  + \omega_*  ^2 }}e^{\beta_*  t} \big[\beta_*  \cos \omega_*  t +
 \omega_*  \sin \omega_*  t\big]\\
&\quad  - ik \times \frac{{c_{12} (k)}}{{\beta_*  ^2  + \omega_*  ^2
}} e^{\beta_*  t} \left[\beta_*  \sin \omega_*  t - \omega_*  \cos
\omega_*  t\right]
 \end{split}
\end{equation}
with
\begin{equation}\label{4.41}
\begin{split}
[&{c_{10}},~   {c_{11}},~   {c_{12}}]^T =\frac{1} {{3{\sigma_* ^2} +
2\sigma_*  + 1 +   {{\left| k \right|}^2}}}\\
&\left[ {\begin{array}{*{20}{c}}
   {\sigma_*I_3} & {\sigma_*(\sigma_*+1)I_3}&{(\sigma_*+1)ik\times} \\
   {-\sigma_*I_3} & {(2\sigma_*^2+\sigma_*+\left| k \right|^2+1)I_3 } & { -(\sigma_*+1)ik\times}  \\
   {\frac{\frac{3}{2}\sigma_*^2+\frac{3}{2}\sigma_*+1+\left| k \right|^2}{\omega_*}I_3 } &
{\frac{(\sigma_*+1)(\sigma_*+1+\left| k \right|^2)}{2\omega_*}I_3 }
& {
\frac{\frac{3}{2}\sigma_*^2+\frac{1}{2}+\left| k \right|^2}{\omega_*} ik\times}  \\
 \end{array} } \right]\left[ {\begin{array}{*{20}{c}}
   {{{M}_{1,0}}}  \\
   { {{M}_{2,0}}}  \\
   {{{M}_{3,0}}}  \\
 \end{array} } \right].
 \end{split}\end{equation}

Now, one can obtain the explicit representation of $\hat{U}$
$=[\hat{\rho}$, $~\hat{u}$,~$\hat{\Theta}$,~$\hat{E}$,~$\hat{B}]$ in
the following from the previous computations.
\begin{theorem}\label{thm4.3}
Suppose $U=[\rho,~u,~\Theta,~E,~B]$ to be the solution of the
initial problem \eqref{4.1}-\eqref{4.2} on the linearized
homogeneous equations with initial data $U_0$
$=[\rho_{0}$,~$u_{0}$,~$\Theta_0$,~$E_0$,~$B_0]$ which satisfies
\eqref{4.3}. For $(t,k)\in(0,\infty)\times\mathbb{R}^3 $ with
$|k|\neq0$, we have the decomposition
\begin{equation}\label{4.42}
\left[ {\begin{array}{*{20}c}
   {\hat{\rho }  (t,k)}  \\
   {\hat{u}  (t,k)}  \\
   {\hat{\Theta} (t,k)}  \\
   {\hat{E}(t,k)}  \\
   {\hat{B}(t,k)}  \\
\end{array}} \right] = \left[ {\begin{array}{*{20}c}
   {\hat{\rho}  (t,k)}  \\
   {\hat{u}_{ ||} (t,k)}  \\
    {\hat{\Theta}  (t,k)}  \\
   {\hat{E}_{||} (t,k)}  \\
   0  \\
\end{array}} \right] + \left[ {\begin{array}{*{20}c}
   0  \\
   {\hat{u}_{  \bot } (t,k)}  \\
          0\\
   {\hat{E}_ \bot  (t,k)}  \\
   {\hat{B}_ \bot  (t,k)}  \\
\end{array}} \right],
\end{equation}
here $\hat{u}_{ ||},\hat{u}_{  \bot }$ are defined as
$$\hat{u}_{
||}=\tilde{k}\tilde{k}\cdot\hat{u},\ ~ \hat{u}_{ \bot
}=-\tilde{k}\times(\tilde{k}\times\hat{u})=(I_3-\tilde{k}\otimes\tilde{k})
\hat{u},$$ and likewise for $\hat{E}_{||},\hat{E}_ \bot$ and
$\hat{B}_ \bot.$ Denote
\begin{equation}\label{4.43}
\left[ {\begin{array}{*{20}c}
   {M_{1 } (t,k)}  \\
   {M_2 (t,k)}  \\
   {M_3 (t,k)}  \\
\end{array}} \right]: = \left[ {\begin{array}{*{20}c}
   {\hat{u}_{  \bot } (t,k)}  \\
   {\hat{E}_ \bot  (t,k)}  \\
   {\hat{B}_ \bot  (t,k)}  \\
\end{array}} \right],{\rm{  }}\left[ {\begin{array}{*{20}c}
   {M_{1 ,0} (k)}  \\
   {M_{2,0} (k)}  \\
   {M_{3,0} (k)}  \\
\end{array}} \right]: = \left[ {\begin{array}{*{20}c}
   {\hat{u}_{ 0, \bot } (t,k)}  \\
   {\hat{E}_{0, \bot } (t,k)}  \\
   {\hat{B}_{0, \bot } (t,k)}  \\
\end{array}} \right].
\end{equation}
 Then, there exit matrices $G^I_{8\times8}(t,k)$ and
$G^{II}_{9\times9}(t,k)$ such that
\begin{equation}\label{4.44}
\left[ {\begin{array}{*{20}c}
   {\hat{\rho} (t,k)}  \\
      {\hat{u}_{||} (t,k)}  \\
{\hat{\Theta} (t,k)}  \\
   {\hat{E}_{||} (t,k)}  \\
\end{array}} \right] = G_{8 \times 8}^I (t,k)\left[ {\begin{array}{*{20}c}
   {\hat{\rho} _{0} (k)}  \\
      {\hat{u}_{0,||} (k)}  \\
{\hat{\Theta} _{0} (k)}  \\
   {\hat{E}_{0,||} (k)}  \\
\end{array}} \right]
\end{equation}
and
\begin{equation}\label{4.45}
\left[ {\begin{array}{*{20}c}
   {M_{1} (t,k)}  \\
      {M_2 (t,k)}  \\
   {M_3 (t,k)}  \\
\end{array}} \right] = G_{9 \times 9}^{II} (t,k)\left[ {\begin{array}{*{20}c}
   {M_{1,0} (k)}  \\
     {M_{2,0} (k)}  \\
   {M_{3,0} (k)}  \\
\end{array}} \right],
\end{equation}
where $ G_{8 \times 8}^I$ is explicitly determined by
representations \eqref{4.23}, \eqref{4.33}, \eqref{4.29},
\eqref{4.24} for $\hat{\rho}(t,k)$, $\hat{u_{||}(t,k)}$,
$\hat{\Theta}(t,k)$, $\hat{E_{||}(t,k)}$ with $c_i(k)$, $(1\leq
i\leq 9)$ defined by \eqref{4.26}, \eqref{4.30}, \eqref{4.34} in
terms of $\hat{\rho}_0(k)$, $\hat{u}_{||,0}(k)$,
$\hat{\Theta}_0(k)$, $\hat{E}_{||,0}(k)$ $
\left(\tilde{k}i|k|^{-1}\hat{\rho}_0\right)$; and $G_{9 \times
9}^{II} $ is determined by the representations \eqref{4.38},
\eqref{4.39} \eqref{4.40} for $M_{1}(t,k),$ $M_{2}(t,k),$
$M_{3}(t,k)$ with $c_{10}(k)$ , $c_{11}(k)$ and $c_{12}(k)$ defined
by \eqref{4.41} in terms of $M_{1,0}(k),$ $M_{2,0}(k),$
$M_{3,0}(k)$.
\end{theorem}
\vspace{0.3cm} \noindent 4.4 \textbf{Refined $L^p-L^q$ time-decay
property.} We utilize Theorem \ref{thm4.3} to acquire some refined
$L^p-L^q$ time-decay property for every component of the solution
$U$ $=[\rho$,~$u$,~$\Theta$,~$E$,~$B]$ in this subsection. For this
purpose, we first search the subtle time-frequency pointwise
estimates on $\hat{U}$
$=[\hat{\rho}$,~$\hat{u}$,~$\hat{\Theta}$,~$\hat{E}$,~$\hat{B}]$ as
follows
\begin{lemma}\label{L4.3}
Suppose $U=[\rho,~u,~\Theta,~E,~B]$ to be the solution of the
linearized homogeneous equations \eqref{4.1} with initial data $U_0$
$=[ \rho _{0}$, $u _{0}$, $\Theta_0$, $E _0$, ${B}_0]$ which satisfy
\eqref{4.3}. Then, there exist constants $\gamma>0,C>0$ such that
for $(t,k)\in(0,\infty)\times\mathbb{R}^3 $,
\begin{equation}\label{4.46}
|\hat{\rho}   (t,k)| \leq C
e^{-\frac{t}{2}}\left|[\hat{\rho}_0,~\hat{u}_0,~\hat{\Theta}_0]\right|,
\end{equation}
\begin{equation}\label{4.47}
\begin{split}
  \left| {\hat u(t,k)} \right| &\leqslant C{e^{ - \frac{t}
{2}}}\left| {\left[ {{{\hat \rho }_0}(t,k),{{\hat u}_0}(t,k),{{\hat
\Theta }_0}(t,k),{{\hat E}_0}(t,k)}
 \right]} \right|  \\
  & + C\left| {\left[ {{{\hat u}_0}(t,k),{{\hat E}_0}(t,k),{{\hat B}_0}(t,k)} \right]} \right| \cdot
    \left\{ {\begin{array}{*{20}{c}}
   {{e^{ - \gamma t}} + \left| k \right|{e^{ - \gamma {{\left| k \right|}^2}t}}}
    & {if~\left| k \right| \leqslant 1,}  \\
   {{e^{ - \gamma t}} + \frac{1}
{{\left| k \right|}}{e^{\frac{{ - \gamma t}}
{{{{\left| k \right|}^2}}}}}} & {if~\left| k \right| > 1,}  \\
 \end{array} } \right.
\end{split}
\end{equation}
\begin{equation}\label{4.48}
|\hat{\Theta}   (t,k)| \leq C
e^{-\frac{t}{2}}\left|[\hat{\rho}_0,~\hat{u}_0,~\hat{\Theta}_0]\right|,
\end{equation}
\begin{equation}\label{4.49}
\begin{split}
  \left| {\hat E(t,k)} \right| &\leqslant C{e^{ - \frac{t}
{2}}}\left| {\left[ {{{\hat u}_0}(t,k),{{\hat \Theta
}_0}(t,k),{{\hat E}_0}(t,k)}
 \right]} \right|  \\
  & + C\left| {\left[ {{{\hat u}_0}(t,k),{{\hat E}_0}(t,k),{{\hat B}_0}(t,k)} \right]} \right| \cdot
    \left\{ {\begin{array}{*{20}{c}}
   {{e^{ - \gamma t}} + \left| k \right|{e^{ - \gamma {{\left| k \right|}^2}t}}}
    & {if~\left| k \right| \leqslant 1,}  \\
   {{\left| k \right|^{-2}e^{ - \gamma t}} +
{e^{\frac{{ - \gamma t}}
{{{{\left| k \right|}^2}}}}}} & {if~\left| k \right| > 1,}  \\
 \end{array} } \right.
\end{split}
\end{equation}
and
\begin{equation}\label{4.50}
\begin{split}
  \left| {\hat B(t,k)} \right| &\leq
   C\left| {\left[ {{{\hat u}_0}(t,k),{{\hat E}_0}(t,k),{{\hat B}_0}(t,k)} \right]} \right| \cdot
    \left\{ {\begin{array}{*{20}{c}}
   {{\left| k \right|e^{ - \gamma t}} + {e^{ - \gamma {{\left| k \right|}^2}t}}}
    & {if~\left| k \right| \leqslant 1,}  \\
   {{\left| k \right|^{-1}e^{ - \gamma t}} +
{e^{\frac{{ - \gamma t}}
{{{{\left| k \right|}^2}}}}}} & {if~\left| k \right| > 1,}  \\
 \end{array} } \right.
\end{split}
\end{equation}
\end{lemma}
\noindent \emph{Proof.} Firstly, we search the upper bound of
$\hat{\rho}$ defined as \eqref{4.23}.
 In fact, from Lemma \ref{L4.1}, it is directly to check
 \eqref{4.26} to obtain
\[\left[ {\begin{array}{*{20}{c}}
   {{c_1}}  \\
   {{c_2}}  \\
   {{c_3}}  \\
 \end{array} } \right] = \left[ {\begin{array}{*{20}{c}}
   {O(1){{\left| k \right|}^2}} & { - O(1){{\left| k \right|}^3}i} & { - O(1){{\left| k \right|}^2}}  \\
   {O(1)} & {O(1){{\left| k \right|}^3}i} & {O(1){{\left| k \right|}^2}}  \\
   {O(1)} & { - O(1)\left| k \right|} & { - O(1){{\left| k \right|}^2}}  \\
 \end{array} } \right]\left[ {\begin{array}{*{20}{c}}
   {{{\hat \rho }_0}}  \\
   {\tilde k \cdot {{\hat u}_0}}  \\
   {{{\hat \Theta }_0}}  \\
 \end{array} } \right]\]
as $|k|\rightarrow0$, and
\[\left[ {\begin{array}{*{20}{c}}
   {{c_1}}  \\
   {{c_2}}  \\
   {{c_3}}  \\
 \end{array} } \right] = \left[ {\begin{array}{*{20}{c}}
   {O(1)} & { - O(1){{\left| k \right|}^{ - 1}}i} & { - O(1)}  \\
   {O(1)} & {O(1){{\left| k \right|}^{ - 1}}i} & {O(1)}  \\
   {O(1){{\left| k \right|}^{ - 1}}} & { - O(1)i} & {O(1){{\left| k \right|}^{ - 1}}}  \\
 \end{array} } \right]\left[ {\begin{array}{*{20}{c}}
   {{{\hat \rho }_0}}  \\
   {\tilde k \cdot {{\hat u}_0}}  \\
   {{{\hat \Theta }_0}}  \\
 \end{array} } \right]\]
as $|k|\rightarrow\infty$. Then, after putting the previous
computations into \eqref{4.26}, one has
\[\begin{split}
  \hat \rho (t,k) =& \left( {O(1){{\left| k \right|}^2}{{\hat \rho }_0} - O(1)
  {{\left| k \right|}^3}i\tilde k \cdot {{\hat u}_0} - O(1){{\left| k \right|}^2}
  {{\hat \Theta }_0}} \right){e^{\sigma t}}   \\
   &+ \left( {O(1){{\hat \rho }_0} + O(1){{\left| k \right|}^3}i\tilde k \cdot {{\hat u}_0}
   + O(1){{\left| k \right|}^2}{{\hat \Theta }_0}} \right){e^{\beta~ t}}\cos \omega t  \\
   &+ \left( {O(1){{\hat \rho }_0} - O(1)\left| k \right|i\tilde k \cdot {{\hat u}_0} - O(1)
   {{\left| k \right|}^2}{{\hat \Theta }_0}} \right){e^{\beta~ t}}\sin \omega t,  \\
\end{split} \]
as $|k|\rightarrow0$, and
\[\begin{split}
  \hat \rho (t,k) =& \left( {O(1){{\hat \rho }_0} - O(1){{\left| k \right|}^{ - 1}}i\tilde k \cdot
  {{\hat u}_0} - O(1){{\hat \Theta }_0}} \right){e^{\sigma t}}   \\
   &+ \left( {O(1){{\hat \rho }_0} + O(1){{\left| k \right|}^{ - 1}}i\tilde k \cdot {{\hat u}_0}
   + O(1){{\hat \Theta }_0}} \right){e^{\beta~ t}}\cos \omega t  \\
   &+ \left( {O(1){{\left| k \right|}^{ - 1}}{{\hat \rho }_0} - O(1)i\tilde k \cdot {{\hat u}_0}
   + O(1){{\left| k \right|}^{ - 1}}{{\hat \Theta }_0}} \right){e^{\beta~ t}}\sin \omega t, \\
\end{split} \]
as $|k|\rightarrow\infty$. Therefore, one can obtain \eqref{4.46}.
Similarly, one can get \eqref{4.48} and the first term on the right
hand side of \eqref{4.49}.

 Next, we search the upper bound of
$\hat{u}_{||}(t,k)$ defined as \eqref{4.33}. In fact, from Lemma
\ref{L4.1}, it is directly to check \eqref{4.34} to get
\[\left[ {\begin{array}{*{20}{c}}
   {{c_7}}  \\
   {{c_8}}  \\
   {{c_9}}  \\
 \end{array} } \right] = \left[ {\begin{array}{*{20}{c}}
   {O(1)} & { - O(1)} & {O(1)\left| k \right|i} & { - O(1){{\left| k \right|}^2}}  \\
   { - O(1)} & {O(1)} & { - O(1)\left| k \right|i} & {O(1){{\left| k \right|}^2}}  \\
   {O(1)} & { - O(1)} & { - O(1)i} & {O(1)\left( {1 - {{\left| k \right|}^2}} \right)}  \\
 \end{array} } \right]\left[ {\begin{array}{*{20}{c}}
   {{{\hat \rho }_0}}  \\
   {\tilde k \cdot {{\hat u}_0}}  \\
   {{{\hat \Theta }_0}}  \\
   {\tilde k \cdot {{\hat E}_0}}  \\
 \end{array} } \right]\]
as $|k|\rightarrow0$, and
\[\left[ {\begin{array}{*{20}{c}}
   {{c_7}}  \\
   {{c_8}}  \\
   {{c_9}}  \\
 \end{array} } \right] = \left[ {\begin{array}{*{20}{c}}
   {O(1)\left( {1 - {{\left| k \right|}^{ - 1}}i} \right)} & { - O(1)} & {O(1){{\left| k \right|}^{ - 1}}i}  \\
   {O(1)\left( {{{\left| k \right|}^{ - 2}} - {{\left| k \right|}^{ - 1}}i} \right)} & {O(1)} & { - O(1){{\left| k \right|}^{ - 1}}i}  \\
   {O(1)\left( {{{\left| k \right|}^{ - 1}} - i} \right)} & { - O(1){{\left| k \right|}^{ - 1}}} & { - O(1){{\left| k \right|}^{ - 1}}i}  \\
 \end{array} } \right]\left[ {\begin{array}{*{20}{c}}
   {{{\hat \rho }_0}}  \\
   {\tilde k \cdot {{\hat u}_0}}  \\
   {{{\hat \Theta }_0}}  \\
 \end{array} } \right]\]
as $|k|\rightarrow\infty$. Therefore, after putting the previous
computations into \eqref{4.33}, one has
\[\begin{split}
  \tilde k \cdot \hat u(t,k) = &\left( {O(1){{\hat \rho }_0} - O(1)\tilde k \cdot {{\hat u}_0} + O(1)
  \left| k \right|i{{\hat \Theta }_0} - O(1){{\left| k \right|}^2}i\tilde k \cdot {{\hat E}_0}}
  \right){e^{\sigma t}}  \\
   &+ \left( { - O(1){{\hat \rho }_0} + O(1)\tilde k \cdot {{\hat u}_0} - O(1)\left| k \right|
   i{{\hat \Theta }_0} + O(1){{\left| k \right|}^2}i\tilde k \cdot {{\hat E}_0}} \right){e^{\beta~ t}}
   \cos \omega t  \\
   &+ \left( {O(1){{\hat \rho }_0} - O(1)\tilde k \cdot {{\hat u}_0} - O(1)i{{\hat \Theta }_0} + O(1)
   \left( {1 - {{\left| k \right|}^2}} \right)\tilde k \cdot {{\hat E}_0}} \right){e^{\beta~ t}}\sin \omega t,  \\
\end{split} \]
as $|k|\rightarrow0$, and
\[\begin{split}
  \tilde k \cdot \hat u(t,k) =& \left( {O(1)\left( {1 - {{\left| k \right|}^{ - 1}}i} \right){{\hat \rho }_0}
  - O(1)\tilde k \cdot {{\hat u}_0} - O(1){{\left| k \right|}^{ - 1}}i{{\hat \Theta }_0}} \right){e^{\sigma t}} \\
   &+ \left( {O(1)\left( {{{\left| k \right|}^{ - 2}} + {{\left| k \right|}^{ - 1}}i} \right){{\hat \rho }_0}
   + O(1)\tilde k \cdot {{\hat u}_0} - O(1){{\left| k \right|}^{ - 1}}i{{\hat \Theta }_0}} \right){e^{\beta~ t}}
   \cos \omega t \\
   &+ \left( {O(1)\left( {{{\left| k \right|}^{ - 1}} - i} \right){{\hat \rho }_0} - O(1){{\left|
    k \right|}^{ - 1}}\tilde k \cdot {{\hat u}_0} - O(1){{\left| k \right|}^{ - 1}}i{{\hat \Theta }_0}}
     \right){e^{\beta~ t}}\sin \omega t,
\end{split} \]
as $|k|\rightarrow\infty$. Then, from above computations, one obtain
the first term on the right hand side of \eqref{4.47}. Similarly, we
get \eqref{4.50} and the second term on the right hand side of both
\eqref{4.47} and \eqref{4.49}. Now, we have finished the proof of
Lemma \ref{L4.3}. \hfill $\Box$ \vspace{0.2cm}

With the help of Lemma \ref{L4.3}, one can refine the time-decay
property for the solution $U=$ $[\rho$, $u$, $\Theta$, $E$, $B]$
obtained in Theorem \ref{thm4.2} in the following.
\begin{theorem}\label{thm4.4}
Let $1\leq p,r\leq 2\leq q\leq\infty,l\geq0$ and  an integer
$m\geq0$. Assume $U(t)=e^{tL}U_0$ is the solution of the initial
problem \eqref{4.1}-\eqref{4.2} with initial data $[\rho_0$, $u_0$,
$\Theta_0$, $E_0$, $B_0]$ which satisfies \eqref{4.3}. Then,
$U=[\rho,~u,~\Theta,~E,~B]$ satisfies
\begin{equation}\label{4.51}
{\left\| {{\nabla ^m}\rho \left( t \right)} \right\|_{{L^q}}}
\leqslant C{e^{ - \frac{t} {2}}}\left( {{{\left\| {\left[ {{\rho
_0},{u_0},{\Theta _0}} \right]} \right\|}_{{L^p}}} + {{\left\|
{{\nabla ^{m + {{\left[ {3\left( {\frac{1} {r} - \frac{1} {q}}
\right)} \right]}_ + }}}\left[ {{\rho _0},{u_0},{\Theta _0}}
\right]} \right\|}_{{L^r}}}} \right),
\end{equation}

\begin{equation}\label{4.52}
\begin{split}
  {\left\| {{\nabla ^m}u\left( t \right)} \right\|_{{L^q}}} \leqslant& C{e^{ - \frac{t}
{2}}}\left( {{{\left\| {\left[ {{\rho _0},{\Theta _0}} \right]}
\right\|}_{{L^p}}} + {{\left\| {{\nabla ^{m + {{\left[ {3\left(
{\frac{1} {r} - \frac{1}
{q}} \right)} \right]}_ + }}}\left[ {{\rho _0},{\Theta _0}} \right]} \right\|}_{{L^r}}}} \right)  \\
  & + C{\left( {1 + t} \right)^{ - \frac{3}
{2}\left( {\frac{1} {p} - \frac{1} {q}} \right) - \frac{{m + 1}}
{2}}}{\left\| {\left[ {{u_0},{E_0},{B_0}} \right]} \right\|_{{L^p}}}  \\
  & + C{\left( {1 + t} \right)^{ - \frac{{l + 1}}
{2}}}{\left\| {{\nabla ^{m + {{\left[ {l + 3\left( {\frac{1} {r} -
\frac{1}
{q}} \right)} \right]}_ + }}}\left[ {{u_0},{E_0},{B_0}} \right]} \right\|_{{L^r}}},  \\
\end{split}
\end{equation}
\begin{equation}\label{4.53}
{\left\| {{\nabla ^m}\Theta \left( t \right)} \right\|_{{L^q}}}
\leqslant C{e^{ - \frac{t} {2}}}\left( {{{\left\| {\left[ {{\rho
_0},{u_0},{\Theta _0}} \right]} \right\|}_{{L^p}}} + {{\left\|
{{\nabla ^{m + {{\left[ {3\left( {\frac{1} {r} - \frac{1} {q}}
\right)} \right]}_ + }}}\left[ {{\rho _0},{u_0},{\Theta _0}}
\right]} \right\|}_{{L^r}}}} \right),
\end{equation}
\begin{equation}\label{4.54}
\begin{split}
\left\|\nabla^mE(t)\right\|_{L^q}\leq
&C(1+t)^{-\frac{3}{2}(\frac{1}{p}-\frac{1}{q})-\frac{m+1}{2}}\left\|{[u_{0},\Theta_0,E_0,B_0]  }\right\|_{L^p}\\
&+C(1+t)^{-\frac{l}{2}}\left\|{\nabla^{m+[l + 3(\frac{1}{r} -
\frac{1}{q})]_ +}[u_{0},\Theta_0,E_0,B_0] }\right\|_{L^r},
\end{split}
\end{equation}
and
\begin{equation}\label{4.55}
\begin{split}
\left\|\nabla^mB(t)\right\|_{L^q}\leq
&C(1+t)^{-\frac{3}{2}(\frac{1}{p}-\frac{1}{q})-\frac{m}{2}}\left\|{[u_{0},E_0,B_0]  }\right\|_{L^p}\\
&+C(1+t)^{-\frac{l}{2}}\left\|{\nabla^{m+[l + 3(\frac{1}{r} -
\frac{1}{q})]_ +}[u_{0},E_0,B_0] }\right\|_{L^r},
\end{split}
\end{equation}
for $t\geq0$, here $C=C(p,q,r,l,m)$ and $[l + 3(\frac{1}{r} -
\frac{1}{q})]_ + $ is defined as \eqref{4.15}.
\end{theorem}
\noindent \emph{Proof.}~We only give the estimate for $\Theta$. Take
$1\leq p,r\leq 2\leq q\leq\infty$ and an integer $m\geq0$, it
follows from Hausdorff-Young inequality with $\frac{1} {{q'}}$ $+
\frac{1} {q} =1$ and \eqref{4.48} that
\begin{equation}\label{C1}
\begin{split}{\left\| {{\nabla ^m}\Theta (t)} \right\|_{L_x^q}}& \leq C
{\left\| {{|k| ^m}\Theta (t)} \right\|_{L_k^{q'}}}\\
& \leqslant C{e^{ - \frac{t} {2}}}\left( {{{\left\| {{{\left| k
\right|}^m}\left[ {{{\hat \rho }_0},{{\hat u}_0},{{\hat \Theta }_0}}
\right]} \right\|}_{{L^{q'}}\left( {|k| \leqslant 1} \right)}} +
{{\left\| {{{\left| k \right|}^m}\left[ {{{\hat \rho }_0},{{\hat
u}_0},{{\hat \Theta }_0}} \right]} \right\|}_{{L^{q'}}\left( {|k|
\geqslant 1} \right)}}} \right). \end{split}\end{equation}
Now let us estimate each term on the right hand side of \eqref{C1}.
For the first term, fixing $\varepsilon>0$ sufficiently small and
using Holder inequality $\frac{1} {{q'}} = \frac{1} {{p'}} +
\frac{{p' - q'}} {{p'q'}}$ with  $\frac{1} {{p'}}$ $+ \frac{1} {p}
=1$, one has
\begin{equation}\notag
\begin{split}
  {\left\| {{{\left| k \right|}^m}\left[ {{{\hat \rho }_0},{{\hat u}_0},
  {{\hat \Theta }_0}} \right]} \right\|_{{L^{q'}}\left( {|k| \leq 1} \right)}} &=
   {\left\| {{{\left| k \right|}^{ - \frac{{p' - q'}}
{{p'q'}}\left( {3 - \varepsilon } \right)}}{{\left| k \right|}^{m +
\frac{{p' - q'}} {{p'q'}}\left( {3 - \varepsilon } \right)}}\left[
{{{\hat \rho }_0},{{\hat u}_0},{{\hat \Theta }_0}} \right]}
 \right\|_{{L^{q'}}\left( {|k| \leq 1} \right)}}  \\
 &  \leqslant \left\| {{{\left| k \right|}^{ - \left( {3 - \varepsilon } \right)}}}
  \right\|_{{L^1}\left( {|k| \leq 1} \right)}^{\frac{{p' - q'}}
{{p'q'}}}{\left\| {{{\left| k \right|}^{m + \frac{{p' - q'}}
{{p'q'}}\left( {3 - \varepsilon } \right)}}\left[ {{{\hat \rho
}_0},{{\hat u}_0},{{\hat \Theta }_0}} \right]}
 \right\|_{{L^{p'}}\left( {|k| \leq 1} \right)}}  \\
  & \leqslant C{\left\| {{{\left| k \right|}^{m + \frac{{p' - q'}}
{{p'q'}}\left( {3 - \varepsilon } \right)}}\left[ {{{\hat \rho
}_0},{{\hat u}_0},{{\hat \Theta }_0}}
 \right]} \right\|_{{L^{p'}}\left( {|k| \leq 1} \right)}}  \\
  & \leq C {{\left\| {\left[ {{\hat{\rho}
_0},{\hat{u}_0},{\hat{\Theta} _0}} \right]} \right\|}_{{L^{p'}}\left( {|k| \leq 1} \right) }}\\
& \leq C {{\left\| {\left[ {{{\rho} _0},{{u}_0},{{\Theta} _0}}
\right]} \right\|}_{{L^{p}}}}.
\end{split} \end{equation}
For the second term, using Holder inequality $\frac{1} {{q'}} =
\frac{1} {{r'}} + \frac{{r' - q'}} {{r'q'}}$ with $\frac{1} {{r'}}$
$+ \frac{1} {r} =1$ and $\varepsilon>0$ small enough, we have
\begin{equation}\notag
\begin{split}
  {\left\| {{{\left| k \right|}^m}\left[ {{{\hat \rho }_0},{{\hat u}_0},
  {{\hat \Theta }_0}} \right]} \right\|_{{L^{q'}}\left( {|k| \geqslant 1} \right)}} &=
   {\left\| {{{\left| k \right|}^{ - \frac{{r' - q'}}
{{r'q'}}\left( {3 + \varepsilon } \right)}}{{\left| k \right|}^{m +
\frac{{r' - q'}}
{{r'q'}}\left( {3 + \varepsilon } \right)}}\left[ {{{\hat \rho }_0},{{\hat u}_0},{{\hat \Theta }_0}} \right]}
 \right\|_{{L^{q'}}\left( {|k| \geqslant 1} \right)}}  \\
 &  \leqslant \left\| {{{\left| k \right|}^{ - \left( {3 + \varepsilon } \right)}}}
  \right\|_{{L^1}\left( {|k| \geqslant 1} \right)}^{\frac{{r' - q'}}
{{r'q'}}}{\left\| {{{\left| k \right|}^{m + \frac{{r' - q'}}
{{r'q'}}\left( {3 + \varepsilon } \right)}}\left[ {{{\hat \rho }_0},{{\hat u}_0},{{\hat \Theta }_0}} \right]}
 \right\|_{{L^{r'}}\left( {|k| \geqslant 1} \right)}}  \\
  & \leqslant C{\left\| {{{\left| k \right|}^{m + \frac{{r' - q'}}
{{r'q'}}\left( {3 + \varepsilon } \right)}}\left[ {{{\hat \rho }_0},{{\hat u}_0},{{\hat \Theta }_0}}
 \right]} \right\|_{{L^{r'}}\left( {|k| \geqslant 1} \right)}}  \\
  & \leqslant C{\left\| {{{\left| k \right|}^{m + {{\left[ {3\left( {\frac{1}
{r} - \frac{1}
{q}} \right)} \right]}_ - } + 1}}\left[ {{{\hat \rho }_0},{{\hat u}_0},{{\hat \Theta }_0}} \right]}
 \right\|_{{L^{r'}}\left( {|k| \geqslant 1} \right)}}\\
 & \leq C {{\left\|
{{\nabla ^{m + {{\left[ {3\left( {\frac{1} {r} - \frac{1} {q}}
\right)} \right]}_ + }}}\left[ {{\rho _0},{u_0},{\Theta _0}}
\right]} \right\|}_{{L^r}}},
\end{split} \end{equation}
which together the above estimate implies \eqref{4.53}. Similarly as
before, we can get \eqref{4.51}, \eqref{4.52}, \eqref{4.54} and
\eqref{4.55}. Then, we finished the proof of Theorem \ref{thm4.4}. \hfill $\Box$\\

Based on Theorem \ref{thm4.4}, we list some particular cases as
follows for later use.
\begin{coro}\label{Corollary4.2}
Assusme $U(t)=e^{tL}U_0$ is the solution of the initial problem
\eqref{4.1}-\eqref{4.2} with initial data $U$ $=[\rho_0$, $u_0$,
$\Theta_0$, $E_0$, $B_0]$ which satisfies \eqref{4.3}. Then,
$U=[\rho,~u,~\Theta,~E,~B]$ satisfies
\begin{equation}\label{4.56}
\left\{
\begin{split}
  &\left\| {\rho \left( t \right)} \right\| \leqslant C{e^{ - \frac{t}
{2}}}\left\| {\left[ {{\rho _0},{u_0},{\Theta _0}} \right]} \right\|,   \\
 & \left\| {u\left( t \right)} \right\| \leqslant C{e^{ - \frac{t}
{2}}}\left\| {\left[ {{\rho _0},{\Theta _0}} \right]} \right\| +
C{\left( {1 + t} \right)^{ - \frac{5}
{4}}}{\left\| {\left[ {{u_0},{E_0},{B_0}} \right]} \right\|_{{L^1} \cap {{\dot H}^2}}},   \\
 & \left\| {\Theta \left( t \right)} \right\| \leqslant C{e^{ - \frac{t}
{2}}}\left\| {\left[ {{\rho _0},{u_0},{\Theta _0}} \right]} \right\|,   \\
 & \left\| {E\left( t \right)} \right\| \leqslant C{\left( {1 + t} \right)^{ - \frac{5}
{4}}}{\left\| {\left[ {{u_0},{\Theta_0},{E_0},{B_0}} \right]} \right\|_{{L^1} \cap {{\dot H}^3}}},   \\
 & \left\| {B\left( t \right)} \right\| \leqslant C{\left( {1 + t} \right)^{ - \frac{3}
{4}}}{\left\| {\left[ {{u_0},{E_0},{B_0}} \right]} \right\|_{{L^1}
\cap {{\dot H}^2}}},
\end{split}\right.
\end{equation}
\begin{equation}\label{4.57}
\left\{\begin{split}& {\left\| {\rho \left( t \right)}
\right\|_{{L^\infty }}} \leqslant C{e^{ - \frac{t} {2}}}{\left\|
{\left[ {{\rho _0},{u_0},{\Theta _0}} \right]} \right\|_{{L^2} \cap
{{\dot H}^2}}},\\
 &  {\left\| {u\left( t \right)}
\right\|_{{L^\infty }}} \leqslant C{e^{ - \frac{t} {2}}}{\left\|
{\left[ {{\rho _0},{\Theta _0}} \right]} \right\|_{{L^1} \cap {{\dot
H}^2}}} + C{\left( {1 + t} \right)^{ - 2}}{\left\| {\left[
{{u_0},{E_0},{B_0}} \right]} \right\|_{{L^1}
\cap {{\dot H}^5}}},  \\
&   {\left\| {\Theta \left( t \right)} \right\|_{{L^\infty }}}
\leqslant C{e^{ - \frac{t}
{2}}}{\left\| {\left[ {{\rho _0},{u_0},{\Theta _0}} \right]} \right\|_{{L^2} \cap {{\dot H}^2}}},  \\
&   {\left\| {E\left( t \right)} \right\|_{{L^\infty }}} \leqslant
C{\left( {1 + t} \right)^{ - 2}}
  {\left\| {\left[ {{u_0},{\Theta_0},{E_0},{B_0}} \right]} \right\|_{{L^1} \cap {{\dot H}^6}}}, \\
&   {\left\| {B\left( t \right)} \right\|_{{L^\infty }}} \leqslant
C{\left( {1 + t} \right)^{ - \frac{3} {2}}}{\left\| {\left[
{{u_0},{E_0},{B_0}} \right]} \right\|_{{L^1} \cap {{\dot H}^5}}},
\end{split}\right.
\end{equation}
and
\begin{equation}\label{4.58}
\left\{\begin{split}
 & \left\| {\nabla B\left( t \right)} \right\| \leqslant C{\left( {1 + t} \right)^{ - \frac{5}
{4}}}{\left\| {\left[ {{u_0},{E_0},{B_0}} \right]} \right\|_{{L^1} \cap {{\dot H}^4}}},  \\
 & \left\| {{\nabla ^s}\left[ {E\left( t \right),B\left( t \right)} \right]} \right\| \leqslant
  C{\left( {1 + t} \right)^{ - \frac{5}
{4}}}{\left\| {\left[ {{u_0},\Theta_0,{E_0},{B_0}} \right]}
\right\|_{{L^2} \cap {{\dot H}^{s+ 3}}}}.
\end{split}\right.
\end{equation}
holds for any $t\geq0$.
\end{coro}
%
\section{ Time-decay rates for system \eqref{2.2}}
In this section, let us give the proof of Proposition \ref{prop2.2}
and Proposition \ref{prop2.3}. For the solution
$U=[\rho,~u,~\Theta,~E,~B]$ of the nonlinear initial problem
\eqref{2.2}-\eqref{2.3}, we search the time-decay rates of the
energy $\|U(t)\|^2_N$ and the high-order energy $ \| \nabla U(t)
\|^2_{N-1}$ in the first two subsections. In the last subsection,
the time-decay rates in $L^q$ with $2\leq q \leq \infty$ for every
component $\rho$, $u$, $\Theta$, $E$ and $B$ of the solution $U$ are
presented.

In the following, since we shall utilize the linear $L^p-L^q$
time-decay property of the homogeneous equations \eqref{4.1}
investigated in the section above to the nonlinear equations
\eqref{2.2}, we rewrite \eqref{2.2} in the following form:
\begin{equation}
\label{5.1}
\left\{\begin{aligned}
&\partial_t \rho + \nabla\cdot u=g_1, \\
&\partial_t u  +\nabla \rho +\nabla \Theta+E+ u   =g_2, \\
&\partial_t \Theta+\frac{2}{3}\nabla\cdot u +  \Theta= g_3 , \\
&  \partial_t E-\nabla\times B-u =g_4, \\
&  \partial_t B+\nabla\times E=0,\\
& \nabla\cdot E=-\rho,\quad \nabla\cdot B=0,\quad
(t,x)\in(0,\infty)\times\mathbb{R}^3,
\end{aligned} \right.
\end{equation}
with
\begin{equation}
 \left\{
 \begin{split}
   &{{g_1} =  - \nabla  \cdot \left( {\rho u} \right),}  \\
 &  {{g_2} =  - \left( {u \cdot \nabla } \right)u - \left( {\frac{{1 + \Theta }}
{{1 + \rho }} - 1} \right)\nabla \rho  - u \times B,}  \\
&   {{g_3} =  - u \cdot \nabla \Theta  - \frac{2}{3}\Theta \nabla
\cdot u + \frac{1}
{3}{{\left| u \right|}^2},}  \\
&   {{g_4} = \rho u.}
 \end{split}
  \right.
\end{equation}
Then, by the Duhamel principle, the solution $U$ can be  formally
written as
\begin{equation}\label{5.3}
U(t) = e^{tL} U_0  + \int_0^t {e^{(t - y)L} [g_1 (y),g_2 (y),g_3
(y),g_4(y),0]} dy,
\end{equation}
here, $ e^{tL}$ is defined as \eqref{4.14}.
\begin{remark}
 In the time integral term of \eqref{5.3}, since $[g_1 (y),g_2 (y),g_3 (y),g_4(y),0]$ satisfies
compatible condition \eqref{4.3}, it makes sense that $e^{(t - y)L}$
acts on $[g_1 (y),g_2 (y),g_3 (y),g_4(y),0]$ for $0\leq y \leq t$.
\end{remark}
\vspace{2mm} \noindent 5.1. \textbf{Decay rate for the energy
functional.} In this subsection, let us search the time-decay
estimate \eqref{2.12} in Proposition \ref{prop2.2} for the energy
$\left\|{U(t)}\right\|^2_N$. We begin with the following Lemma which
can be seen straightforward from the proof of Theorem \ref{thm3.1}.
\begin{lemma}\label{L5.1}
 Suppose $U$ $=[\rho$, $u$, $\Theta$, $E$, $B]$ to be the solution of the
 initial problem \eqref{2.2}-\eqref{2.3} with
 $U_0$ $=[\rho_0$, $u_0$, $\Theta_0$, $E_0$, $B_0]$ which satisfies \eqref{2.4} obtained by Proposition \ref{prop2.1}.
 If $\mathcal {E}_s(U_0)$
is small enough, then for any $t \geq 0$, it holds that
\begin{equation}\label{5.4}
\frac{d}{dt}\mathcal {E}_s(U(t))+\gamma \mathcal {D}_s(U(t))\leq 0.
\end{equation}
\end{lemma}
Based on Lemma \ref{L5.1}, one can check that
\begin{equation}\notag
\begin{split}
(1+t)^l\mathcal {E}_s(U(t))+& \gamma\int_0^t(1+y)^l\mathcal
{D}_s(U(y))dy\\
&\leq \mathcal {E}_{s}(U_0)+l\int_0^t(1+y)^{l-1}\mathcal
{E}_s(U(y)) dy\\
&\leq \mathcal {E}_{s}(U_0)+C l\int_0^t(1+y)^{l-1}\left(\left\|{B(y)
}\right\|^2+\mathcal {D}_{s+1}(U(y)) \right)dy,
\end{split}
\end{equation}
where we have used $\mathcal {E}_s(U(t))\leq \left\|{B(t)
}\right\|^2+\mathcal {D}_{s+1}(U(t))$. Using \eqref{5.4} again, one
has
 $$\mathcal {E}_{s+2}(U(t))+\gamma\int_0^t \mathcal {D}_{s+2}(U(s))dy\leq \mathcal {E}_{s+2}(U_0)$$
and
\begin{equation}\notag
\begin{split}
(1+t)^{l-1}\mathcal {E}_{s+1}(U(t))+&
\gamma\int_0^t(1+y)^{l-1}\mathcal
{D}_{s+1}(U(y))dy\\
&\leq \mathcal
{E}_{s+1}(U_0)+C(l-1)\int_0^t(1+y)^{l-2}\left(\left\|{B(y)
}\right\|^2+\mathcal {D}_{s+2}(U(y)) \right)dy.
\end{split}
\end{equation}
Therefore, by iterating the previous estimates, we have
\begin{equation}\label{5.5}
\begin{split}
(1+t)^l\mathcal {E}_s(U(t))+& \gamma\int_0^t(1+y)^l\mathcal
{D}_s(U(y))dy\\
&\leq C\mathcal {E}_{s+2}(U_0)+C\int_0^t(1+y)^{l-1}\left\|{B(y)
}\right\|^2dy
\end{split}
\end{equation}
for $1<l<2.$

 Now, let us estimate the integral term on the right
hand side of \eqref{5.5}.
Applying the last linear estimate on $B$ in \eqref{4.56} to
\eqref{5.3}, one has \begin{equation}\label{B}
 \left\| {B\left( t \right)}
\right\| \leqslant C {\left( {1 + t} \right)^{ - \frac{3}
{4}}}{\left\| {\left[ {{u_0},{E_0},{B_0}} \right]} \right\|_{{L^1}
\cap {{\dot H}^2}}} + C\int_0^t {{{\left( {1 + t - y} \right)}^{ -
\frac{3} {4}}}{{\left\| {\left[ {{g_2}(y),{g_4}(y)} \right]}
\right\|}_{{L^1} \cap {{\dot H}^2}}}} dy.
\end{equation}
 It is directly to check that for any $0\leq y \leq t$,
 \[{\left\| {\left[ {{g_2}(y),{g_4}(y)} \right]} \right\|_{{L^1} \cap {{\dot H}^2}}}
 \leqslant C\left\| {U(y)} \right\|_3^2 \leqslant C{\mathcal {E}_s}\left( {U(y)} \right) \leqslant
  C{\left( {1 + y} \right)^{ - \frac{3}
{2}}}{\mathcal {E} _{s,\infty }}\left( {U(t)} \right),\]
 where ${\mathcal {E} _{s,\infty }}\left(
{U(t)} \right): = \mathop {\sup }\limits_{0 \leqslant y \leqslant t}
{\left( {1 + y} \right)^{\frac{3} {2}}}{\mathcal {E}_s}\left( {U(y)}
\right).$ Plugging this into \eqref{B} implies
\begin{equation}\label{5.6}
\left\| {B\left( t \right)} \right\| \leqslant C {\left( {1 + t}
\right)^{ - \frac{3} {4}}}\left( {{{\left\| {\left[
{{u_0},{E_0},{B_0}} \right]} \right\|}_{{L^1} \cap {{\dot H}^2}}} +
{\mathcal {E} _{s,\infty }}\left( {U(t)} \right) }
\right).\end{equation}

 Next, we prove the uniform-in-time bound of
${\mathcal {E} _{N,\infty }}\left( {U(t)} \right)$ which implies the
decay rates of the energy functional ${\mathcal {E} _{s }}\left(
{U(t)} \right)$ and thus $\|U(t)\|_s^2$. In fact, by choosing
$l=\frac{3}{2}+\varepsilon$ in \eqref{5.5} with $\varepsilon>0$
sufficiently small and using \eqref{5.6}, it follows that
\begin{equation}\notag\begin{split}
  {\left( {1 + t} \right)^{\frac{3}
{2} + \varepsilon }}{\mathcal {E} _s}\left( {U(t)} \right) &+ \gamma
{\int_0^t {\left( {1 + y} \right)} ^{\frac{3}
{2} + \varepsilon }}{\mathcal {D}_s}\left( {U(y)} \right)dy  \\
 &  \leqslant C{\mathcal {E} _{s + 2}}\left( {{U_0}} \right) +
C{\left( {1 + t} \right)^\varepsilon }\left( {\left\| {\left[
{{u_0},{E_0},{B_0}} \right]} \right\|_{_{{L^1} \cap {{\dot H}^2}}}^2
+ {{\left[ {{\mathcal {E} _{s,\infty }}\left( {U(t)} \right)}
\right]}^2}} \right),
\end{split} \end{equation}
which yields
$${\left(
{1 + t} \right)^{\frac{3} {2}}}{\mathcal {E}_s}\left( {U(t)} \right)
\leqslant C\left( {{\mathcal {E}_{s + 2}}\left( {{U_0}} \right) +
\left\| {\left[ {{u_0},{E_0},{B_0}} \right]} \right\|_{{L^1}}^2 +
{{\left[ {{\mathcal {E} _{s,\infty }}\left( {U(t)} \right)}
\right]}^2}} \right),
$$
 and
thus
$$
{\mathcal {E} _{s,\infty }}\left( {U(t)} \right) \leqslant C\left(
{{\epsilon _{s + 2}}
  {{\left( {{U_0}} \right)}^2} + {{\left[ {{\mathcal {E}_{s,\infty }}\left( {U(t)} \right)}
  \right]}^2}} \right), $$
since $\epsilon _{s + 2}\left( {{U_0}} \right)>0 $ is small enough,
it holds that ${\mathcal {E} _{s,\infty }}\left( {U(t)} \right)
\leqslant C{\epsilon _{s + 2}}{\left( {{U_0}} \right)^2} $ for any
$t\geq 0,$ which gives ${\left\| {U(t)} \right\|_s} \leqslant
C{\mathcal {E}_s}{\left( {U(t)} \right)^{\frac{1} {2}}} \leqslant
C{\epsilon _{s + 2}}\left( {{U_0}} \right){\left( {1 + t} \right)^{
- \frac{3} {4}}}$, that is \eqref{2.12}. \hfill $\Box$
 \vspace{0.3cm}

  \noindent5.2. \textbf{Decay rate for the high-order energy
functional.} In this subsection, we search the decay estimate of the
high-order energy functional $\|\nabla U(t)\|^2_{s-1}$, that is
\eqref{2.13} in Proposition \ref{prop2.2}. For that, we are reduced
to establish the time-decay estimates on $\|\nabla B\|$ and
$\|\nabla^s[E,B]\|$ with the help of the following Lemma.
\begin{lemma}\label{L5.3}
Suppose $U=[\rho,~u,~\Theta,~E,~B]$ to be the solution of the
initial problem \eqref{2.2}-\eqref{2.3} with $U_0$ $=[\rho_0$,
$u_0$, $\Theta_0$, $E_0$, $B_0]$ which satisfies \eqref{2.4}
obtained in Proposition \ref{prop2.1}. If $\mathcal {E}_s(U_0)$ is
small enough, then there exit the high-order energy functional
$\mathcal {E}_s^h(\cdot) $ and the high-order dissipation rate
$\mathcal {D}_s^h(\cdot) $ such that
\begin{equation}\label{5.7}
\frac{d}{dt}\mathcal {E}_s^h(U(t))+\gamma \mathcal {D}_s^h(U(t))\leq
C\|\nabla B\|^2,
\end{equation}
holds for any $t\geq 0.$
\end{lemma}
\noindent \emph{Proof.}  The proof is very similar to the proof of
Theorem \ref{thm3.1}. In fact, after letting $|\alpha| \geq 1$, then
corresponding to \eqref{3.3}, \eqref{3.14}, \eqref{3.26} and
\eqref{3.27},  it can also be tested that
\begin{equation}\notag
\frac{d}{dt}\sum_{1\leq|\alpha|\leq s} [\langle A^{I
}_0(W_I)\partial^\alpha W_{I},\partial^\alpha
    W_{I}\rangle+\|\partial^\alpha W_{II}\|^2
    ]+\|\nabla u\|_{s-1}^2+\frac{1}{3}\|\nabla\Theta\|_{s-1}^2\leq
    C\|W\|_s\|W_I\|_s^2,\ \ \
\end{equation}
\begin{equation}\notag
 \begin{split}  \frac{d}{dt}\sum_{1\leq\alpha\leq s-1}\langle \frac{1}{2(1+\rho)}\partial ^\alpha\rho & -\partial
   ^\alpha\nabla\cdot u,
  \partial ^\alpha \rho\rangle+
    \gamma\|\nabla\rho\|_{s-1}^2\\
  & \leq C ( \|W\|_s  \|W_I\|_s^2+  \|\nabla^2 u\|_{s-2}^2+\|\nabla\Theta\|_{s-1}^2
  ),
        \end{split}
\end{equation}
\begin{equation}\notag
 \begin{split}  &\frac{d}{dt}\sum_{1\leq|\alpha|\leq s-1}\langle \partial ^\alpha u,\partial ^\alpha E\rangle  +
    \gamma\|\nabla E\|_{s-2}^2 \\
    & \leq  C\left(\|\nabla u\|_{s-1}^2 +\|\nabla^2\Theta\|_{s-2}^2\right)
      +C  \|U\|_s\|U\|_s^2+C \|\nabla^2 u\|_{s-2}\|\nabla
      B\|_{s-2}
        \end{split}
\end{equation}
and
\begin{equation}\notag
\frac{d} {{dt}}\sum\limits_{1\leq\left| \alpha  \right| \leqslant s-
2} {\left\langle { - \nabla  \times {\partial ^\alpha }E,{\partial
^\alpha }B} \right\rangle }  + \gamma \left\| {\nabla^2 B}
\right\|_{s - 3}^2 \leqslant C\left\| {[\nabla u, \nabla^2 E]}
\right\|_{s - 3}^2 + C\left\| U  \right\|_{_s}^2\left\|
{\nabla[\rho, u]} \right\|_{_{s - 1}}^2.\end{equation}

Now, in the similar way as in \emph{Step 5} of Theorem \ref{thm3.1}.
We define the high-order energy functional as
\begin{equation}\label{00}
 \begin{split}
 \mathcal{E}_s^h(U(t))=& \sum_{1\leq|\alpha|\leq s} [\langle A^{I
}_0(W_I)\partial^\alpha W_{I},\partial^\alpha
    W_{I}\rangle+\|\partial^\alpha W_{II}\|^2
    ]\\
    & +\mathcal {K}_1 \sum\limits_{1\leq\left| \alpha  \right| \leqslant s -
1} \langle \frac{1}{2(1+\rho)}\partial ^\alpha\rho-\partial
^\alpha\nabla\cdot u,
  \partial ^\alpha \rho\rangle
\\&+\mathcal {K}_2 \sum\limits_{1\leq\left| \alpha  \right| \leqslant
s - 1}\langle \partial ^\alpha u,\partial ^\alpha E\rangle +\mathcal
{K}_3 \sum\limits_{1\leq\left| \alpha  \right| \leqslant s - 2}
{\left\langle { - \nabla  \times {\partial ^\alpha }E,{\partial
^\alpha }B} \right\rangle }.
  \end{split}
\end{equation}
Similarly, one can take $0<\mathcal {K}_3\ll\mathcal
{K}_2\ll\mathcal {K}_1\ll 1$ be sufficiently small with $\mathcal
{K}_2^{\frac{3}{2}}\ll\mathcal {K}_3$, such that
$\mathcal{E}_s^h(U(t))\sim \|\nabla U(t)\|_{s-1}^2 $, that is,
$\mathcal{E}_s^h(\cdot)$ is really a high-order energy functional
which satisfies \eqref{2.6}, and moreover, the summation of the
 four above estimates with coefficients corresponding to
\eqref{00} implies \eqref{5.7} with $\mathcal {D}^h_s(\cdot)$
defined as \eqref{2.8}. We have finished the
proof of Lemma \ref{L5.3}. \hfill $\Box$\\

Based on Lemma \ref{L5.3}, one can check that
\begin{equation}\notag
\frac{d}{dt}\mathcal {E}_s^h(U(t))+\gamma \mathcal {E}_s^h(U(t))\leq
C\left(\|\nabla B\|^2+\|\nabla^s[E,B]\|^2\right),
\end{equation}
which implies
\begin{equation}\label{5.8}
\begin{split}
\mathcal {E}_s^h(U(t))\leq  \mathcal {E}_s^h(U_0)e^{-\gamma t}
+ C\int_0^t{e^{-\gamma (t-y)}\left(\|\nabla B(y)\|^2+\|\nabla^s
[E(y),~B(y)] \|^2 \right)}dy.
\end{split}
\end{equation}

 Now, let us estimate the time integral term on the right hand side of the
previous inequality. Firstly, we have
\begin{lemma}\label{L5.4}
Suppose $U=[\rho,~u,~\Theta,~E,~B]$ to be the solution of the
initial problem \eqref{2.2}-\eqref{2.3} with $U_0$ $=[\rho_0$,
$u_0$, $\Theta_0$, $E_0$, $B_0]$ which satisfies \eqref{2.4}
obtained in Proposition \ref{prop2.1}. If $\epsilon_{N+6}(U_0)$ is
small enough, then
\begin{equation}\label{5.9}
{\left\| {\nabla B(t)} \right\|^2} + {\left\| {{\nabla ^s}\left[
{E(t),B(t)} \right]} \right\|^2} \leqslant C{\epsilon _{s +
6}}{({U_0})^2}{(1 + t)^{ - \frac{5} {2}}},
\end{equation}
for any $t\geq 0.$
\end{lemma}
\noindent \emph{Proof.} The proof is similar to that of the
isentropic case in \cite{Duan11}. Apply the first linear estimate on
$\nabla B$ in \eqref{4.58} to \eqref{5.3} so that
\begin{equation}\notag
\begin{split}
  \left\| {\nabla B\left( t \right)} \right\| & \leqslant C{\left( {1 + t} \right)^{ - \frac{5}
{4}}}{\left\| {\left[ {{u_0},{E_0},{B_0}} \right]} \right\|_{{L^1}
\cap {{\dot H}^4}}}+ C\int_0^t{\left( {1 +t- y} \right)^{ - \frac{5}
{4}}}{\left\| {\left[ {{g_2(y)},{g_4(y)}} \right]} \right\|_{{L^1} \cap {{\dot H}^4}}}dy \\
 & \leqslant C{\left( {1 + t} \right)^{ - \frac{5}
{4}}}{\left\| {\left[ {{u_0},{E_0},{B_0}} \right]} \right\|_{{L^1}
\cap {{\dot H}^4}}}+ C\int_0^t{\left( {1 +t- y} \right)^{ - \frac{5}
{4}}}{\left\| {U(y)} \right\|^2_{\max\{5,s\}  }}dy\\
 & \leqslant C{\left( {1 + t} \right)^{ - \frac{5}
{4}}}{\left\| {\left[ {{u_0},{E_0},{B_0}} \right]} \right\|_{{L^1}
\cap {{\dot H}^4}}}+ C\int_0^t{\left( {1 +t- y} \right)^{ - \frac{5}
{4}}}\epsilon_{s+6}(U_0)^2   \left( {1 + y} \right)^{ - \frac{3}
{2}} dy\\
 & \leqslant C\epsilon_{s+6}(U_0)   \left( {1 + t} \right)^{ - \frac{5}
{4}}.
\end{split}
\end{equation}
Similarly, applying the second linear estimate on $ {{\nabla
^s}\left[ {E\left( t \right),B\left( t \right)} \right]} $ in
\eqref{4.58} to \eqref{5.3}, one has
\begin{equation}\notag
\begin{split}
  & \left\| {{\nabla ^s}\left[ {E\left( t \right),B\left( t \right)} \right]}
  \right\|\\
  & \leqslant
  C{\left( {1 + t} \right)^{ - \frac{5}
{4}}}{\left\| {\left[ {{u_0},\Theta_0,{E_0},{B_0}} \right]}
\right\|_{{L^2} \cap {{\dot H}^{s + 3}}}}+C\int_0^t{\left( {1 + t-y}
\right)^{ - \frac{5} {4}}}{\left\| {\left[
{{g_2(y)},g_3(y),{g_4(y)}} \right]} \right\|_{{L^2} \cap {{\dot
H}^{s + 3}}}}dy\\
  & \leqslant
  C{\left( {1 + t} \right)^{ - \frac{5}
{4}}}{\left\| {\left[ {{u_0},\Theta_0,{E_0},{B_0}} \right]}
\right\|_{{L^2} \cap {{\dot H}^{s + 3}}}}+C\int_0^t{\left( {1 + t-y}
\right)^{ - \frac{5} {4}}}{\left\| { {U(y)} } \right\|^2_{  {{s +
4}}}}dy\\
  & \leqslant
  C{\left( {1 + t} \right)^{ - \frac{5}
{4}}}{\left\| {\left[ {{u_0},\Theta_0,{E_0},{B_0}} \right]}
\right\|_{{L^2} \cap {{\dot H}^{s + 3}}}}+C\int_0^t{\left( {1 + t-y}
\right)^{ - \frac{5} {4}}}\epsilon_{s+6}(U_0)^2   \left( {1 + y}
\right)^{ - \frac{3}
{2}} dy\\
 & \leqslant C\epsilon_{s+6}(U_0)   \left( {1 + t} \right)^{ - \frac{5}
{4}}.
\end{split}
\end{equation}

Where we have used \eqref{2.12} and the smallness of
$\epsilon_{s+6}(U_0)  $. Now, we have finished the proof of Lemma
\ref{L5.4}. \hfill $\Box$\\

 Then, by putting \eqref{5.9} into \eqref{5.8},
we have
\begin{equation}\notag
\mathcal {E}_N^h(U(t))\leq \mathcal {E}_N^h(U_0)e^{-\lambda t}+ C
\epsilon_{N+6}(U_0)^2(1+t)^{-\frac{5}{2}}.
\end{equation}
Since $\mathcal {E}^h_N(U(t))\sim \|\nabla U(t) \|^2_{N-1}$ holds
true for any $t\geq 0$, \eqref{2.16} follows. Therefore, we have
finished the proof of Proposition \ref{2.2}. \vspace{0.3cm}

 \noindent
 5.3 \textbf{Decay rate in $L^q$ .}  In this subsection, we will search the decay rates of solutions
$U$ $=[\rho$, $u$, $\Theta$, $E$, $B]$ in $L^q$ with $2\leq q\leq
+\infty$ to the initial problem \eqref{2.2}-\eqref{2.3} by proving
Proposition \ref{prop2.3}. Throughout this subsection, we always
assume that $\epsilon_{13}(U_0)>0$ is small enough. First, for
$s\geq 4$, Proposition \ref{prop2.2} shows that if
$\epsilon_{s+2}(U_0)$ is small enough,
\begin{equation}\label{5.10}
\|U(t)\|_s\leq C \epsilon_{s+2}(U_0)(1+t)^{-\frac{3}{4}},
\end{equation}
and if $\epsilon_{s+6}(U_0)$ is small enough,
\begin{equation}\label{5.11}
\|\nabla U(t)\|_{s-1}\leq C \epsilon_{s+6}(U_0)(1+t)^{-\frac{5}{4}}.
\end{equation}
Now, let us establish the estimates on $B$, $[u,~E]$ and
$[\rho,~\Theta]$ in the following.\\
 \noindent \emph{Estimate on $\|B\|_{L^q}$.} For $L^2$ rate, it is
 straightforward
 from \eqref{5.10} to obtain
 $$\|B(t)\|\leq C \epsilon_{~6} (U_0)(1+t)^{-\frac{3}{4}}.$$
 For $L^\infty$ rate, by utilizing $L^\infty$ estimate on $B$
 in \eqref{4.57} to \eqref{5.3}, we have
$${\left\| {B(t)} \right\|_{{L^\infty }}} \leqslant C{(1 + t)^{ -
\frac{3} {2}}}{\left\| {[{u_0},{E_0},{B_0}]} \right\|_{{L^1} \cap
{{\dot H}^5}}} + C\int_0^t {{{(1 + t - y)}^{ - \frac{3} {2}}}}
{\left\| {[{g_2}(y),{g_4}(y)]} \right\|_{{L^1} \cap {{\dot
H}^5}}}dy.$$
 Because of \eqref{5.10},
$${\left\| {[{g_2}(t),{g_4}(t)]} \right\|_{{L^1} \cap {{\dot H}^5}}}
\leqslant C\left\|
 {U(t)} \right\|_6^2 \leqslant C{\epsilon _8}{({U_0})^2}{(1 + t)^{ - \frac{3}
{2}}},$$
 one has
 $${\left\| {B(t)} \right\|_{{L^\infty }}} \leqslant C\epsilon_8 ({U_0}) {(1 + t)^{ -
 \frac{3}
{2}}}.$$ Therefore, by $L^2-L^\infty$ interpolation
\begin{equation}\label{5.12}
\|B(t)\|_{L^q}\leq C \epsilon_8 ({U_0}) {(1 + t)^{ - \frac{3}
{2}+\frac{3}{2q}}},
\end{equation}
for $2\leq q \leq \infty.$\\
\noindent\emph{Estimate on $\|[u,E]\|_{L^q}$.} For $L^2$ rate,
utilizing the $L^2$ estimate on $u$ and $E$ in \eqref{4.56} to
\eqref{5.3}, we have
\begin{equation}\notag
\begin{split}
  \left\| {u\left( t \right)} \right\| \leqslant& C{\left( {1 + t} \right)^{ - \frac{5}
{4}}}\left( {\left\| {\left[ {{\rho _0},{\Theta _0}} \right]}
\right\| + {{\left\|
{\left[ {{u_0},{E_0},{B_0}} \right]} \right\|}_{{L^1} \cap {{\dot H}^2}}}} \right)  \\
   &+ C\int_0^t {{{\left( {1 + t - y} \right)}^{ - \frac{5}
{4}}}\left( {\left\| {\left[ {{g_1}(y),{g_3}(y)} \right]} \right\| +
{{\left\| {\left[ {{g_2}(y),{g_4}(y)} \right]} \right\|}_{{L^1} \cap
{{\dot H}^2}}}} \right)dy}
\end{split}
\end{equation}
and $$\left\| {E\left( t \right)} \right\| \leqslant C{\left( {1 +
t} \right)^{ - \frac{5} {4}}}{\left\| {\left[
{{u_0},\Theta_0,{E_0},{B_0}} \right]} \right\|_{{L^1} \cap {{\dot
H}^3}}} + C\int_0^t {{{\left( {1 + t - y} \right)}^{ - \frac{5}
{4}}}{{\left\| {\left[ {{g_2}(y),{g_3}(y),{g_4}(y)} \right]}
\right\|}_{{L^1} \cap {{\dot H}^3}}}dy.} $$ Because of \eqref{5.10},
$$\left\| {\left[ {{g_1}(t),{g_3}(t)} \right]} \right\| + {\left\|
{\left[ {{g_2}(t),{g_3}(t),{g_4}(t)} \right]} \right\|_{{L^1} \cap
{{\dot H}^3}}} \leqslant C\left\| {U(t)} \right\|_4^2 \leqslant
C{\epsilon _6}{({U_0})^2}{(1 + t)^{ - \frac{3} {2}}},$$ it holds
that
\[\left\| {\left[ {u\left( t \right),E\left( t \right)} \right]} \right\| \leqslant C{\epsilon _6}({U_0}){\left( {1 + t} \right)^{ - \frac{5}
{4}}}.\] For $L^\infty$ rate, utilizing the $L^\infty$ estimates on
$u$ and $E$ in \eqref{4.57} to \eqref{5.3}, we have
\begin{equation}\notag
\begin{split}
  {\left\| {u\left( t \right)} \right\|_{{L^\infty }}} \leqslant& C{\left( {1 + t} \right)^{ - 2}}
  \left( {{{\left\| {\left[ {{\rho _0},{\Theta _0}} \right]} \right\|}_{{L^1} \cap {{\dot H}^2}}}
   + {{\left\| {\left[ {{u_0},{E_0},{B_0}} \right]} \right\|}_{{L^1} \cap {{\dot H}^5}}}} \right)  \\
   &+ C\int_0^t {{{\left( {1 + t - y} \right)}^{ - 2}}\left( {{{\left\| {\left[ {{g_1}(y),{g_3}(y)}
   \right]} \right\|}_{{L^1} \cap {{\dot H}^2}}} + {{\left\| {\left[ {{g_2}(y),{g_4}(y)} \right]}
   \right\|}_{{L^1} \cap {{\dot H}^5}}}} \right)} dy
\end{split} \end{equation}
and $${\left\| {E\left( t \right)} \right\|_{{L^\infty }}} \leqslant
C{\left( {1 + t} \right)^{ - 2}}{\left\| {\left[
{{u_0},{\Theta_0},{E_0},{B_0}} \right]} \right\|_{{L^1} \cap {{\dot
H}^6}}} + C\int_0^t {{{\left( {1 + t - y} \right)}^{ - 2}}{{\left\|
{\left[ {{g_2}(y),{g_3}(y),{g_4}(y)} \right]} \right\|}_{{L^1} \cap
{{\dot H}^6}}}dy}. $$ Since
\begin{equation}\notag
\begin{split}
  {\left\| {\left[ {{g_1}(t),{g_3}(t)} \right]} \right\|_{{{\dot H}^2}}}
   + {\left\| {\left[ {{g_2}(t),{g_3}(t),{g_4}(t)} \right]} \right\|_{{{\dot H}^5} \cap {{\dot H}^6}}}
    \leqslant C\left\| {\nabla U(t)} \right\|_6^2 \leqslant {\epsilon _{13}}{({U_0})^2}{(1 + t)^{ - \frac{5}
{2}}},
\end{split} \end{equation}
and
$$\|\left[ {{g_1}(t),{g_2}(t),{g_3}(t),{g_4}(t)} \right]\|_{L^1}\leq
C \|U(t)\|(\|u(t)\|+\|\nabla U\|)\leq
\epsilon_{10}(U_0)^2(1+t)^{-2},$$ then, it holds that
 $$\|[u(t),~E(t)]\|_{L^\infty}\leq C \epsilon_{13}(U_0)^2(1+t)^{-2}.$$
 Then, by $L^2-L^\infty$ interpolation
\begin{equation}\label{5.13}
\|[u(t),~E(t)]\|_{L^q}\leq C \epsilon_{13} ({U_0}) {(1 + t)^{ -
2+\frac{3}{2q}}},
\end{equation}
for $2\leq q \leq \infty.$\\
\noindent\emph{Estimate on $\|[\rho,\Theta]\|_{L^q}$.} For $L^2$
rate, utilizing the $L^2$ estimates on $\rho$ and $\Theta$ in
\eqref{4.56} to \eqref{5.3}, we have
\begin{equation}\label{5.14}
\begin{split}\left\| {\left[ {\rho ,\Theta } \right]} \right\| \leqslant C{e^{ - \frac{t}
{2}}}\left\| {\left[ {{\rho _0},{u_0},{\Theta _0}} \right]} \right\|
+ C\int_0^t {{e^{ - \frac{{t - y}} {2}}}\left\| {\left[
{{g_1}(y),{g_2}(y),{g_3}(y)} \right]} \right\|dy.} \end{split}
\end{equation}
Since
\begin{equation}\notag
\begin{split}
 \left\| {\left[ {{g_1}(t),{g_2}(t),{g_3}(t)} \right]} \right\| \leqslant C\left\| {\nabla U(t)}
  \right\|_1^2 + \left\| {u(t)} \right\|
  \left({\left\| { {B(t)} } \right\|_{{L^\infty }}}+\|\nabla U(t)\|\right)
   \leqslant C{\epsilon _{10}}{({U_0})^2}{\left( {1 + t} \right)^{ - \frac{5}
{2}}},
 \end{split}
\end{equation}
then \eqref{5.14} implies the slower decay estimate
\begin{equation}\label{5.15}
\begin{split}\left\| {\left[ {\rho(t) ,\Theta(t) } \right]} \right\| \leqslant C
{\epsilon _{~10}}{({U_0})}{\left( {1 + t} \right)^{ - \frac{5}
{2}}}.\end{split}
\end{equation}
Furthermore, after estimating $\left\| {\left[
{{g_1}(t),{g_2}(t),{g_3}(t)} \right]} \right\|$ and utilizing the
previous slower decay estimate, one has
\begin{equation}\notag
\begin{split}
  &\left\| {\left[ {{g_1}(t),{g_2}(t),{g_3}(t)} \right]} \right\|\\
  & \leqslant C{\left\|
   u \right\|_{{L^\infty }}}\left( {\left\| {\nabla \left[ {\rho (t),u(t),\Theta (t)}
    \right]} \right\| + \left\| {\left[ {B(t),u(t)} \right]} \right\|} \right)
   + C\left\| {\left[ {\rho (t),\Theta (t)} \right]} \right\|{\left\| {\nabla
   \left[ {\rho (t),u(t)} \right]} \right\|_2}\\
   & \leqslant C{\epsilon _{13}}{({U_0})^2}{(1 + t)^{ - \frac{{11}}
{4}}}, \end{split}
\end{equation}
it follows from \eqref{5.14} that $$\left\| {\left[ {\rho (t),\Theta
(t)} \right]} \right\| \leqslant C{\epsilon _{13}}{({U_0})}{(1 +
t)^{ - \frac{{11}} {4}}}.$$ For $L^\infty$ rate, by utilizing the
$L^\infty$ estimates on $\rho$ and $\Theta$ in \eqref{4.57} to
\eqref{5.3},
\begin{equation}\label{5.16}
{\left\| {\left[ {\rho ,\Theta } \right]} \right\|_{{L^\infty }}}
 \leqslant  C{e^{ - \frac{t} {2}}}{\left\|
{\left[ {{\rho _0},{u_0},{\Theta _0}} \right]} \right\|_{{L^2} \cap
{{\dot H}^2}}} + C\int_0^t {{e^{ - \frac{{t - y}} {2}}}{{\left\|
{\left[ {{g_1}(y),{g_2}(y),{g_3}(y)} \right]} \right\|}_{{L^2} \cap
{{\dot H}^2}}}dy} .
\end{equation} It is straightforward to check
\begin{equation}\label{5.17}
\begin{split}
  &{\left\| {\left[  {{g_1}(t), {g_2}  (  t),{g_3}(t)} \right]} \right\|_{{L^2}\cap {{\dot H}^2}}}\\
   &\leqslant C{\left\|
  {\nabla U(t)} \right\|_4}\left({\left\| {\left[ {\rho(t),\Theta(t)} \right]} \right\|}+
   {{{\left\| {\left[ {u(t),B(t)} \right]} \right\|}_{{L^\infty }}} + {{\left\|
   {\nabla \left[ {\rho (t),u(t),\Theta (t)} \right]} \right\|}_{{L^\infty }}}} \right)  \\
  &\leq C{\epsilon _{13}}{({U_0})}{(1 +
t)^{ - \frac{{11}} {4}}},\end{split}
\end{equation}
which yields from \eqref{5.16} that
\[{\left\| {\left[ {\rho (t),\Theta (t)} \right]} \right\|_{{L^\infty }}} \leq C{\epsilon _{13}}{({U_0})}{(1 +
t)^{ - \frac{{11}} {4}}}.\]  Then, by $L^2-L^\infty$ interpolation
\begin{equation}\label{5.18}
\|[\rho(t),~\Theta(t)]\|_{L^q}\leq C \epsilon_{13} ({U_0}) {(1 +
t)^{ - \frac{11}{4}}},
\end{equation}
for $2\leq q \leq \infty.$ Therefore, \eqref{5.18},  \eqref{5.13}
and
 \eqref{5.12} give  \eqref{2.14},  \eqref{2.15} and  \eqref{2.16}, respectively. Now, we have finished
  the proof of Proposition \ref{prop2.3}.\hfill $\Box$\\

\noindent {\bf Acknowledgments}
 \vspace{2mm}

This work is supported by China 973 Program(Grant no. 2011CB808002),
the Grants NSFC 11071009 and PHR-IHLB 200906103 and Foundation
Project of Doctor Graduate Student Innovation of Beijing University
of Technology of China.

\vspace{2mm}
 \noindent {\bf Appendix A}
\appendix
\makeatletter
\renewcommand\theequation{A.\@arabic\c@equation }
\makeatother \setcounter{equation}{0}

 \vspace{2mm}

In a medium, the Euler-Maxwell system including natural collision
terms is written as the following (nonconservative)form (see
\cite{Chen84,Dink05,Jero03,Jero05}):
\begin{equation}
{{\partial _t}n + \nabla  \cdot (nu) = 0,}
\end{equation}
\begin{equation}
\label{A.2}
 {m\left[ {{\partial _t}(nu) + \nabla  \cdot (nu \otimes
u)} \right] + k\nabla (n\theta )
 =  - qn(E + u \times B) - \frac{{mnu}}{{{\tau _p}}},}
\end{equation}
\begin{equation}
{{\partial _t}\theta  + u \cdot \nabla \theta  + \frac{2}{3}\theta
\nabla  \cdot u
 = \frac{{{k_0}}}{n}\nabla  \cdot (n\nabla \theta ) - \frac{{2m{{\left| u \right|}^2}}}{{3k}}
 \left( {\frac{1}{{2{\tau _\omega }}} - \frac{1}{{{\tau _p}}}} \right)
 - \frac{1}{{{\tau _\omega }}}(\theta  - {\theta _*}),}
\end{equation}
\begin{equation}
{\epsilon {\partial _t}E - {\mu ^{ - 1}}\nabla  \times B =
\frac{q}{m}nu,} \qquad {\epsilon \nabla  \cdot E = b(x) -
\frac{q}{m}n,}
\end{equation}
\begin{equation}
{{\partial _t}B + \nabla  \times E = 0,}\qquad  {\nabla  \cdot B =
0,}
\end{equation}
for $(t,x)\in (0,+\infty)\times\mathbb{R}^3$. Here,
$\epsilon>0,\mu>0$ and $m>0$ are the permittivity of the medium, the
permeability of the medium and the particle mass, respectively. In
vacuum, $\epsilon=\epsilon_0, \mu=\mu_0$ with
$c=(\epsilon_0\mu_0)^{-\frac{1}{2}}$ being the speed of light. $q$
is the electronic charge, $m$ is the effective electron mass, $k$ is
Boltzmann's constant, $\tau_p$ is the momentum relaxation time,
$\tau_\omega$ is the energy relaxation time, $k_0$ is a constant
multiplier (with the variable density) of heat conduction. The
function $\theta_*(x)$  is the ambient device temperature, and the
function $b(x)$ stands for the prescribed density of positive
charged background ions (doping profile). In this paper, we assume
$q=m=k=\epsilon=\mu=\tau_p=\tau_\omega=\theta_*(x)=b(x)=1$ and
$k_0=0$ for the sake of simplicity. It's well known that \eqref{A.2}
is equivalent to
\begin{equation}
n\partial_t u + n( u\cdot\nabla) u +n\nabla \theta+\theta\nabla n +
nu = -n(E+   u \times B) , \end{equation}
 then, we obtain the desired Euler-Maxwell
system of the form \eqref{1.1}.



\begin{thebibliography}{99}
\bibitem{Chen84}
F. Chen, Introduction to Plasma Physics and Controlled Fusion. Vol.
1, Plenum Press, New York, 1984.

\bibitem{CJW00} G. Q.~Chen, J. W.~Jerome, D. H.~Wang.
Compressible Euler-Maxwell equations, {Transport Theory and
Statistical Physics}, 29 (2000) 311-331.

\bibitem{Dink05}
Andreas Dinklage, et al. Plasma Physics, in: Lect. Notes Phys. Vol.
670, Springer, Berlin, Heidelberg, 2005.

\bibitem{Duan11}  R. J.~Duan. Global smooth flows for the compressible Euler-Maxwell system: Relaxation
case, {Journal of Hyperbolic Differential Equations.} 8 (2011)
375-413.
\bibitem{DH95}
D. Hoff and K. Zumbrun. {Multi-dimensional diffusion waves for the
Navier-Stokes equations of compressible flow}, Indiana Univ. Math.
J, {44} (1995) 603-676.

\bibitem{Jero03}
J. W. Jerome. The Cauchy problem for compressible
hydrodynamic-Maxwell systems: a local theory for smooth solutions,
Differential and Integral Equations 16 (2003) 1345-1368.

\bibitem{Jero05}
J. W. Jerome. Functional Analytic Methods for Evolution Systems, in:
Contemporary Mathematics. Vol. 371, American Mathematical Society,
Providence, 2005, pp. 193-204.

\bibitem{Ka75} T.~Kato.
The Cauchy problem for quasi-linear symmetric hyperbolic systems,
{Arch. Ration. Mech. Anal.} 58 (1975) 181-205.
\bibitem{SK84}
S. Kawashima. {Systems of a hyperbolic-parabolic composite type,
with applications to the equations of magnetohydrodynamics},
Doctoral Thesis, Kyoto University, 1984.

\bibitem{KM81} S.~Klainerman, A.~Majda.
Singular limits of quasilinear hyperbolic systems with large
parameters and the incompressible limit of compressible fluids,
{Comm. Pure Appl. Math.} 34 (1981) 481-524.

\bibitem{Ma84} A.~Majda.
{Compressible Fluid Flow and Systems of Conservation Laws in Several
Space Variables}, Springer-Verlag, New York, 1984.

\bibitem{Ni78} T.~Nishida.
Nonlinear hyperbolic equations and related topics in fluids
dynamics, Publications Math\'ematiques d'Orsay, Universit\'e
Paris-Sud, Orsay, No. 78-02, 1978.

\bibitem{PW08a} Y. J.~Peng, S.~Wang.
Convergence of compressible Euler-Maxwell equations to
incompressible Euler equations, {Comm. Part. Diff. Equations,} 33
(2008) 349-376.

\bibitem{PW08b} Y. J.~Peng, S.~Wang.
Rigorous derivation of incompressible e-MHD equations from
compressible Euler-Maxwell equations, {SIAM J. Math. Anal.} 40
(2008) 540-565.

\bibitem{PWG11} Y. J.~Peng, S.~Wang, Q. L.~Gu.
Relaxation limit and global existence of smooth solutions of
compressible Euler-Maxwell equations, {SIAM J. Math. Anal.} 43
(2011) 944-970.

\bibitem{Stein}
E. M. Stein. {Singular Integrals and Differentiability Properties of
Functions}, Princeton Mathematical Series, No. 30 Princeton
University Press, Princeton, N.J. 1970 xiv+290 pp.

\bibitem{UKW10}
Y. Ueda, S. Wang, S.~Kawashima.{ Large Global existence and
asymptotic decay of solutions to the Euler-Maxwell system}, preprint
2010.

\bibitem{UK11}
Y. Ueda, S.~Kawashima.{Decay property of regularity-loss type for
the Euler-Maxwell system, to appear in Methods and Applications of
Analysis, 2011.}


\bibitem{YW11}
J. W. Yang, S. Wang. The diffusive relaxation limit of
non-isentropic Euler-Maxwell equations for plasmas, {J. Math. Anal.
Appl.} 380 (2011) 343-353.

\end{thebibliography}
\end{document}